\documentclass[11pt]{amsart}
\usepackage{amsmath,amsthm, amscd, amssymb, amsfonts}
\usepackage{color}

\newtheorem{Lem}{Lemma}[section]
\newtheorem{Prop}[Lem]{Proposition}
\newtheorem{Cor}[Lem]{Corollary}

\newtheorem{Theorem}[Lem]{Theorem}
\theoremstyle{definition}
\newtheorem{Def}[Lem]{Definition}
\newtheorem{Rem}[Lem]{Remark}

\newtheorem{Ex}[Lem]{Example}

\usepackage{hhline}

\newcommand{\co}{\operatorname{co} }

\renewcommand{\_}[1]{_{\left( #1 \right)}}

\newcommand{\ot}{{\otimes}}

\newcommand\g{\mathfrak{g}}
\newcommand\Uc{\mathcal{U}}

\newcommand{\X}{{\mathcal X}}

\newcommand{\rt}{\color{black}}

\newcommand{\ku}{\Bbbk}

\newcommand{\Z}{{\mathbb Z}}
\newcommand{\N}{{\mathbb N}}
\newcommand{\I}{{\mathbb I}}

\newcommand{\C}{{\mathcal C}}
\newcommand{\D}{{\mathcal D}}

\newcommand{\ydh}{{}^H_H\mathcal{YD}}

\newcommand{\ydup}{{}^{\Uc(\D', \lambda')}_{\Uc(\D', \lambda')}\mathcal{YD}}

\newcommand{\is}{\I^{\operatorname{s}}}

\newcommand{\Ss}{{\mathcal S}}

\newcommand{\End}{\operatorname{End}}

\newcommand\ad{\operatorname{ad}}
\newcommand\Hom{\operatorname{Hom}}

\newcommand{\NA}{\mathcal{B}}

\newcommand\id{\operatorname{id}}

\def\pf{\begin{proof}}
\def\epf{\end{proof}}

\numberwithin{equation}{section}\theoremstyle{plain}

\newcommand\ord{{\operatorname{ord}}}

\newcommand\Irr{{\operatorname{Irr}}}

\renewcommand\o{\otimes}

\newcommand\YDG{^{\Gamma}_{\Gamma}\mathcal{YD}}
\newcommand\YDL{^{\Lambda}_{\Lambda}\mathcal{YD}}

\newcommand\YDH{^{H}_{H}\mathcal{YD}}

\newcommand\G{\Gamma}

\newcommand\sw[1]{{}_{(#1)}}
\newcommand\swo[1]{{}^{(#1)}}

\newcommand\U{{\bf {U}}}

\newcommand\on{\triangleright}

\begin{document}


\title[Complete reducibility theorems]{Complete reducibility theorems for modules over pointed Hopf algebras}
\author[andruskiewitsch, radford and schneider]{Nicol\'{a}s Andruskiewitsch}
\address{Facultad de Matem\'{a}tica, Astronom\'{i}a y f\'{i}sica\\
Universidad Nacional de C\'{o}rdoba \\ CIEM - CONICET,
(5000) Ciudad Universitaria \\
C\'{o}rdoba \\Argentina} \email{andrus@mate.uncor.edu}
\author[]{David Radford}
\address{Department of Mathematics, Statistics and Computer Science, M/C 249, University of Illinois at Chicago,
851 S. Morgan Street,
Chicago, IL 60607-7045, USA}
\email{radford@uic.edu}
\author[]{Hans-J\"urgen Schneider}
\address{Mathematisches Institut, Universit\"at M\"unchen, Theresienstr. 39,
D-80333 Munich, Germany}
\email{Hans-Juergen.Schneider@mathematik.uni-muenchen.de}
\thanks{Some results of this paper were obtained during a visit of
H.-J. Schneider at the University of C\'ordoba, partially supported
through a grant of CONICET. The work of N. Andruskiewitsch was partially
supported by CONICET, ANPCyT, Mincyt (C\'{o}rdoba) and Secyt (UNC). Research for this paper by
D. Radford was partially supported by NSA Grant H 98230-06-1-10015 and conducted in part during visits {\rt to}
Munich University. He would like to thank the University for its hospitality.}
\dedicatory{Dedicated to Susan Montgomery for her many contributions to Hopf algebras,
ring theory, and service to the mathematical community}
\begin{abstract} We investigate the representation theory of a large class of pointed Hopf algebras, extending results of Lusztig and others.
We classify all simple modules in a suitable category and determine the weight multiplicities; we establish a complete reducibility
theorem in this category.
\end{abstract}

\maketitle

\setcounter{tocdepth}{1}


\section*{Introduction}

{\rt The} main achievements of the {\rt finite-dimensional} representation theory of
finite-dimen\-sional complex semisimple Lie algebras {\rt include}
\begin{enumerate}
  \item[(a)] the parametrization of the simple {\rt finite-dimensional} representations by their highest weights,

  \item[(b)] the complete reducibility of any {\rt finite-dimensional} representation,

  \item[(c)] the determination of the weight multiplicities (and consequently, of the dimension) of
  a {\rt finite-dimensional} highest weight representation.
\end{enumerate}
These classical results of E. Cartan -- part (a) 1913, and of H.
Weyl -- parts (b), (c) 1926, were generalized in many directions.
They hold for Kac-Moody algebras, with symmetrizable generalized
Cartan matrices, in the context of \emph{integrable}
modules from the category $\mathcal O$ instead of
\emph{finite-dimensional} ones, as shown by V. Kac in 1974; see
\cite{K}.

{\rt Drinfeld and Jimbo introduced two quantum versions} of the universal enveloping algebras of a
finite-dimen\-sional simple Lie algebra $\g$, the formal deformation
$U_{\hbar}(\g)$ \cite{Dr0} and the $q$-analogue $U_q(\g)$
\cite{Ji}. The representation theory of the $q$-analogue $U_q(\g)$
was worked out in \cite{L1}, where highest weight modules were
introduced and parts (a) and (c) above were established. Drinfe'ld
observed in \cite{Dr2} that, since $U_{\hbar}(\g)$ and
$U(\g)[[\hbar]]$ are isomorphic algebras by the {\rt White}head Lemma,
their representation theories are equivalent. He also introduced a
quantum Casimir operator to deal with complete reduciblity. Later,
the representation theory of the $q$-analogue $U_q(\g)$, where
$\g$ is a symmetrizable Kac-Moody algebra, was developed in
\cite{L}, where analogues of the highest weight modules from
\cite{K} were introduced and a complete reducibility theorem was
proved \cite[6.2.2]{L} using a quantum version of the Casimir
operator.

Multiparameter quantized enveloping algebras also {\rt have been defined as
formal deformations and as $q$-analogues}, see \cite{Re}
and \cite{OY} respectively. In the context of formal deformations,
this is just a twist in the sense of Drinfeld of $U_{\hbar}(\g)$,
hence their representation theories are equivalent; see also
\cite{LS}. Besides \cite{OY}, multiparameter $q$-analogues of
enveloping algebras were considered in many papers, where the
number of parameters and the group of group-likes is varied; consult \cite{AE}.
{\rt In the last few years, there was a revival
of this question and the representation theories of 2-parameter
deformations were studied; see \cite{BW1,BW2,BGH1,BGH2}, and \cite{HPR} for more general cases and references therein.}

{\rt In \cite{AS-crelle} a} family of pointed Hopf algebras was introduced
 having a Cartan matrix of finite type as one of
the inputs. This family contains the $q$-analogues $U_q(\g)$ and
their multiparametric variants; in fact, they are close to them
but one parameter of deformation for each connected component and
more general linking relations are allowed. Indeed, the family
contains also the parabolic Hopf subalgebras of the $U_q(\g)$'s
and eventually, any pointed Hopf algebra with generic braiding and
finitely generated abelian group of group-likes which is a domain
of finite Gelfand-Kirillov dimension belongs to it \cite{AA,
AS-crelle}. The goal of the present paper is to study the
representation theory of these Hopf algebras {\rt and of their natural generalizations with arbitrary symmetrizable Cartan matrices}. Let us summarize our
main results.

In {\em Sections 1 and 2} we first {\rt consider very general} pointed Hopf algebras
$\Uc(\D,\lambda)$ (see Definition \ref{linkingHopf}) defined {\rt for} a
general YD-datum {\rt $\D$} {\rt over an arbitray abelian group $\G$} and a linking parameter {\rt $\lambda$}.  A YD-datum is a realization of a
diagonal braiding, without further restrictions. We prove a
general structure result, Theorem \ref{Levi}, {\rt analogous }  to the
classical Levi decomposition for Lie algebras.

We focus in {\em{\rt S}ection  \ref{sect:perfect}} on a
special class of linking parameters, the perfect linkings.
Following \cite{RS2}, we introduce the notions of reduced data
$\D_{red}$ and their linking parameters $\ell$. The Hopf algebras
$\Uc(\D,\lambda)$  with perfect linking admit an alternative
presentation in terms of reduced data: $\Uc(\D,\lambda) \simeq
\Uc(\D_{red},\ell)$; this stress{\rt es} the similarity with $U_q(\g)$.
Indeed, $\U  = \Uc(\D_{red},\ell)$ (see Definition \ref{Unondegenerate})  has generators similar to those
of  $U_q(\g)$ (except for the group $\G$ of group-likes that is
more general) and {\rt by Theorem \ref{form}} it is a quotient of a Drinfeld double. {\rt From this description we derive a basic bilinear form $(\;,\;) : \U^- \times \U^+ \to \ku$.  Once the existence of the form is shown,} our approach to establish its main properties is close to previous work in the
literature, specially \cite{L}. We end this section  with a
discussion of data of Cartan type \cite{AS-crelle}.

In {\em{\rt S}ection  \ref{sect:$U$-modules}},  we study the
representation theory of $\U  = \Uc(\D_{red},\ell)$ with
{\rt $\D_{red}$ generic (see Definition \ref{YD-datum}), regular (see Definition \ref{Defregular}), and of Cartan type (see Definition \ref{reducedCartanYD-datum})}.

{\rt The Hopf algebra $\U = \U_L$ in Lusztig's book \cite{L}, where the root datum is $X$-regular, is a special case of our $\U$ above.
The 2-parameter deformations mentioned above and the (generic) multiparameter deformations we have seen in the literature are all special cases of the Hopf algebra $\U$ in this Section.

We extend the results of Chapters 3, 4 and 6 in \cite{L} to our more general context following the strategy of \cite{L}.}
Similarly to \cite{L},
we consider the category $\C$ of $\U$-modules with weight
decomposition (with respect to the action of $\G$), the full
subcategory $\C^{hi}$ (see Subsection \ref{subsection:C}), and the class of integrable modules in $\C$ (see Subsection \ref{subsection:integrable}). {\rt However, Lusztig only considers representations of $\U_L$ with weights of the form $\chi_{\lambda}$, where $\lambda \in X$, and $\chi_{\lambda}(K_{\mu}) = v^{<\mu,\lambda>}$ for all $\mu \in Y$, where the free abelian groups $X,Y$ of finite rank are part of the given root datum in \cite{L}. The weights $\chi_{\lambda}$ (``of type one'') do not make sense in our general context.
Thus our category $\C$ for $\U_L$, where arbitrary elements in $\widehat{\G}$ are allowed as weights, is larger than Lusztig's category $\C$ defined in \cite[3.4.1]{L}.}


{\rt Let $\widehat{\G}^+$ be the set of all dominant characters $\chi \in \widehat{\G}$ for $\D_{red}$ introduced in \cite{RS2}.}
We define Verma modules $M(\chi)$ for all $\chi \in \widehat{\G}$ and Weyl
modules $L(\chi)$ {\rt for all $\chi \in \widehat{\G}^+$.
To define a version of the quantum Casimir operator of \cite[6.11]{L} we have to define a suitable scalar valued function on $\widehat{\G}$ in Lemma \ref{G} extending the integer valued function on $X$ in \cite[Lemma 6.1.5]{L} in the special case considered by Lusztig. It turns out that this is not always possible. But such a function exists if the Dynkin diagram of the generalized Cartan matrix of $\D_{red}$ is connected. We then  reduce the later results to the connected case.

We call an algebra $A$ reductive if all finite-dimensional left $A$-modules are semisimple. If $B \subset A$ is a subalgebra, we say that $A$ is $B$-reductive if all finite-dimensional left $A$-modules which are semisimple over $B$ are semisimple.}
We prove:

\begin{itemize}
    \item [($\alpha$)] {\rt By Corollary \ref{simplebijection} $$\widehat{\G}^+ \to \{ [L] \mid L \in \mathcal{C}^{hi}, L \text{ integrable  and simple} \},\;
\chi \mapsto [L(\chi)],$$ is bijective. Here, $[L]$ denotes the isomorphism class of a module $L$.
This  {\rt establishes}} (a) in the context of integrable modules
in $\C^{hi}$.

\medbreak \item [($\beta$)] {\rt By Theorem \ref{main}} any integrable module in $\C^{hi}$ is completely
reducible. That is, we {\rt prove }(b) in
this context. {\rt This is one of our main results extending Lusztig's Theorem \cite[6.2.2]{L}.
In particular, $\U$ is $\ku\G$-reductive.}

\medbreak \item[($\gamma$)] {\rt Assume that the  braiding is twist-equivalent to a braiding of
Drinfeld-Jimbo type \cite{AS-crelle} (see  \ref{DJ2}). In particular this holds in the case of finite Cartan type.} Then by Theorem \ref{mainsimple}  the weight multiplicities, in particular the
    dimension, of $L(\chi)$ with $\chi$ dominant, are as in the classical
    case. Thus  (c) is {\rt shown}.

\medbreak \item [($\delta$)] {\rt By Theorem
\ref{mainreductive}}
 $\U$ is reductive
iff the index $[\G:\G^2]$
is finite, where $\G^2$ is a subgroup of $\G$ given by the datum $\D_{red}$. Hence {\rt we have determined all Hopf algebras $\U$ satisfying (a) and (b) in the context of finite-dimensional modules.}
\end{itemize}
In the case of 2-parameter deformations of finite Cartan type ($\alpha$) and ($\beta$) have been shown in \cite{BW2} for type A and in \cite{BGH2} for types B, C and D; see also \cite{HPR} for more general cases but under very restrictive assumptions on the braiding.

In {\em Section 5} we {\rt extend } Theorem \ref{mainreductive} {\rt to } the pointed Hopf algebras $\Uc(\D,\lambda)$ {\rt in the case of finite Cartan type. We show in Theorem \ref{pointedreductive} that the linking is perfect if and only if $\Uc(\D,\lambda)$ is $\ku\G$-reductive. Our proof of this close relationship between reductivity and properties of the linking is based on the Levi type decomposition in Theorem \ref{Levi} and recent results on PBW-bases in left coideal subalgebras of quantum groups in \cite{Kh} or \cite{HS}.} Combined with the main results of
\cite{AS-crelle, AA}, {\rt our theory gives  in Theorem \ref{fingrowth-lifting} a characterization of the pointed Hopf algebras $\U$ with finite Cartan matrix and free abelian group of finite rank $\G$ by axiomatic properties.}

{\rt The last author thanks I. Heckenberger for very helpful discussions.}

\section{Nichols algebras and linking}

We denote the ground field by $\ku$, and its multiplicative group of units
by $\ku^{\times}$. {\rt We assume that $\ku$ is algebraically closed of characteristic zero.} By convention, $\mathbb N = \{0, 1, \dots\}$.
If $A$ is an algebra and $g\in A$ is invertible, then $g \on a = gag^{-1}$, $a
\in A$, {\rt denotes the inner automorphism defined by $g$.}

We use standard notation for Hopf algebras; the comultiplication
is denoted $\Delta$ and the antipode $\Ss$. For the first, we use
the Heyneman-Sweedler notation $\Delta(x) = x\_{1}\otimes x\_{2}$. The
left adjoint representation of $H$ on itself is the algebra map
$\ad: H\to \End{\rt (H)}$, $\ad_lx (y) = x\_{1} y \Ss(x\_{2})$, $x,y\in
H$; we shall write $\ad$ for $\ad_l$, omitting  the subscript $l$ unless  strictly needed. There is
also a right adjoint action given by $\ad_r x(y) = \Ss(x\_{1}) y
x\_{2}$. Note that both $\ad_l$ and $\ad_r$ are multiplicative.

\subsection{Yetter-Drinfeld modules and Nichols algebras}

\

For a full account of these structures, the reader is referred to \cite{AS-cambr}.
{\rt Let $H$ be a Hopf algebra with bijective antipode. A
Yetter-Drinfeld module over $H$ is a left $H$-module and a left $H$-comodule with comodule structure denoted by $\delta : V \to H \o V, v \mapsto v\sw{-1} \o v\sw{0},$ such that
$$\delta(h v) = h\sw1 v\sw{-1} \Ss(h\sw3) \o h\sw2 v\sw0$$ for all $v \in V, h \in H$. Let $\ydh$ be the
category of Yetter-Drinfeld modules over $H$ with $H$-linear
and $H$-colinear maps as morphisms.

The category $\ydh$ is monoidal and braided. If $V,W \in {\ydh}$,
then $V \otimes W$ is the tensor product over $\ku$ with the
diagonal action and coaction of $H$ and braiding
$$c_{V,W} : V \otimes W \to W \otimes V, \quad v \otimes w \mapsto v\sw{-1}\cdot w \otimes v\sw0$$
for all $v \in V, w \in W$. This allows us to
consider Hopf algebras in $\ydh$. If $R$ is a Hopf
algebra in the braided category $\ydh$, then the space of primitive elements $P(R)=\{x \in R \mid \Delta(x) = x \o 1 + 1
\o x \}$
is  a Yetter-Drinfeld submodule of $R$.

\smallbreak For $V \in {\ydh}$ the tensor algebra $T(V)=\oplus_{n
\geq 0} T^n(V)$ is an $\N$-graded algebra and coalgebra in the
braided category ${\ydh}$ where the elements of $V=T(V)(1)$ are
primitive. It is a Hopf algebra in ${\ydh}$ since $T(V)(0)=\ku$.}

\smallbreak
We now recall the definition of the fundamental example of {\rt a} Hopf algebra in a category of Yetter-Drinfeld modules.

\begin{Def}\label{DefNichols}
{\rt Let $V \in {\ydh}$ and $I(V) \subset T(V)$  the largest
$\mathbb{N}$-graded ideal and coideal $I \subset T(V)$ such that
$I \cap V =0$. We call $\NA(V) = T(V)/I(V)$ the {\em Nichols
algebra} of $V$. Then $\NA(V) = \oplus_{n \geq 0} \NA^n(V)$ is an
$\mathbb{N}$-graded Hopf algebra in $\ydh$.}
\end{Def}

\begin{Lem}\label{LemmaNichols0}
The Nichols algebra of an object {\rt $V \in {\ydh}$} is (up to isomorphism) the unique $\mathbb{N}$-graded Hopf algebra $R$ in {\rt $\ydh$ }satisfying the following properties:
{\rt \begin{enumerate}
\item $R(0) =\ku,\; R(1) = V$,\label{N1}
\item $R(1)$ generates the algebra $R$,\label{N2}
\item $P(R) = V$.\label{N3}
\end{enumerate}}
If $f : V \to W$ is a morphism in {\rt $\YDH$}, then $T(f)(I(V)) \subset
I(W)$, where $T(f) : T(V) \to T(W)$ is the induced map on the
tensor algebras; hence $f$ induces a morphism between the
corresponding Nichols algebras.
\end{Lem}

\pf \cite[Proposition 2.2, Corollary 2.3]{AS-cambr}. \epf

Let $\G$ be an abelian group. A
Yetter-Drinfeld module over the group algebra {\rt $\ku\G$} is a $\G$-graded vector space
$V = \oplus_{g \in \G} V_g$ which is a left $\G$-module such that
each homogeneous component $V_g$, $g \in \G$, is stable under the
action of $\G$. The $\G$-grading is equivalently described as a
left $\ku\G$-comodule structure $\delta : V \to \ku\G \otimes V$,
$v \mapsto v\sw{-1} \otimes v\sw0$ where $\delta(v) =g \otimes v$
if $v$ is homogeneous of degree $g \in \G$. Let $\YDG$ be the
category of Yetter-Drinfeld modules over $\ku\G$.
{\rt For $V,W \in {\YDG}$ the braiding
is given by $c_{V,W}(v \o w) = g\cdot w \otimes v$
for all $v \in V_g$, $g \in \G$, $w \in W$.}

\subsection{Braided Hopf algebras and
bosonization}\label{subsection:bosonization}

\

{\rt Let $R$ be a
Hopf algebra in $\ydh$.  We denote its comultiplication by $\Delta_R : R \to R \o R, \; r \mapsto r\swo1 \o r\swo2$. The \emph{bosonization} $R \# H$ of $R$ is a Hopf algebra defined as follows. As a vector space, $R\# H = R\otimes
H$; the multiplication and comultiplication of $R\# H$ are given
by the smash-product and smash-coproduct, respectively, that is, for all $r,s \in R, g,h \in H,$
\begin{align*}
(r \o g)(s \o h) &= r (g\sw{1} \cdot s) \o g\sw2h,\\
\delta(r \o g) &= r\swo1 \# (r\swo2){\sw{-1}} g\sw1 \o (r\swo2){\sw0} \# g\sw2.
\end{align*}}

{\rt  Let $\pi: A\to H$ be
a morphism of Hopf algebras. Then $$R = A^{\co H} = \{a\in A \mid
(\id\otimes \pi)\Delta (a) = a\otimes 1\}$$ is the left coideal
subalgebra of  right coinvariants of $\pi$. There is also a left
version: $^{\co H}A =\{a\in A: (\pi\otimes \id)\Delta (a) =
1\otimes a\}$. The subalgebra $A^{\co H} \subset A$ is stable under the left adjoint action $\ad_l$ of $A$, and $^{\co H}A$ is stable under $\ad_r$.

 \begin{Lem}\label{leftright}
Let $\pi: A\to H$ be
a morphism of Hopf algebras, and assume that the antipode $\Ss$ of $A$ is bijective. Then
\begin{enumerate}
\item  $\Ss(A^{\co H}) = {^{\co H}A}$.
\item  Assume that there is a Hopf algebra map $ \iota: H\to A$ such that $\pi\iota = \id_H$. Then the map
$$\kappa : A^{\co H} \to {^{\co H}A},\; a \mapsto \Ss^2(a\sw2) \Ss(\iota\pi(a\sw1)),$$ is bijective, and for all $x,y \in A$,
$\kappa(xy) = \kappa(x\sw2) \ad_r\iota \pi\Ss(x\sw1)(\kappa(y)).$
\end{enumerate}
\end{Lem}
\pf This is easily checked. The map $b \mapsto \Ss^{-2}(b\sw1) \Ss^{-1}(\iota\pi(b\sw2))$ is inverse to $\kappa$.
\epf
\smallbreak Assume the situation of part (2) of Lemma \ref{leftright}. Then $R=A^{\co H}$ is a braided Hopf algebra in $\ydh$,
and the multiplication induces an isomorphism $R\# H \to A, r \# h \mapsto r \iota(h),$ of Hopf algebras.

Conversely any Hopf algebra $R$ in $\ydh$  arises in this way from the bosonization since $\pi= \varepsilon \o \id : R \# H \to H$ is a Hopf algebra map with $(R \# H)^{\co H} = R \o 1$. The braided adjoint action $\ad_c : R \to \End(R)$ is defined for all $x,y \in R$ by
$\ad_c x(y) = \ad x (y)$,
where $\ad$ is the left adjoint action of the bosonization $R \# H$, and where we view $R \to R \# H, \; r \mapsto r \# 1,$ as inclusion.
In particular, if $x \in P(R)$, then $\ad_cx(y) = xy - (x\sw{-1} \cdot y) x\sw0$ is the braided commutator of $x$ and $y$.}

\subsection{Yetter-Drinfeld data}

\

We are interested in {\rt finite-dimensional} Yetter-Drinfeld modules over an abelian group $\G$ that are semisimple as
$\G$-modules; these are described in the following way.

\begin{Def}\label{YD-datum}
A {\em YD-datum}
$\D=\D(\G,(g_i)_{1 \leq i \leq \theta}, (\chi_i)_{1 \leq i \leq \theta})$
consists of an abelian group $\G$, a positive integer $\theta$,
$g_1, \dots,g_{\theta} \in \G$, and characters
$\chi_1,\dots,\chi_{\theta} \in \widehat {\G} = \Hom(\G,\ku
^{\times})$. We let $\I= \{1,2,\dots,\theta\}$ and define
\begin{align}\label{q}
q_{ij}=\chi_j(g_i) \text{ for all } i,j \in \I.
\end{align}
A YD-datum is called {\em generic} if for all $1 \leq i \leq
\theta$, $q_{ii}$ is not a root of unity. We define an equivalence
relation $\sim$ on $\I$, where for all $i,j \in \I$, $i \neq j$,
\begin{align}\label{sim}
i \sim j  \iff \quad
 &\text{There are }i_1,\dots,i_t \in \I, t \geq 2
 \text{ with }i=i_1,j=i_t,\\ &\phantom{a}q_{i_li_{l+1}}q_{i_{l+1}i_l} \neq 1
 \text{ for all }1 \leq l < t.\notag
\end{align}
Let $\X$ be the set of equivalence classes of $\I$ with respect to $\sim$.
\end{Def}

\begin{Rem} The notions of YD-datum and the equivalence relation \eqref{sim} are generalizations
of the notions of Cartan datum and the resulting equivalence relation from \cite[Definition 3.16]{RS2}.

\end{Rem}

Let $\D=\D(\G,(g_i)_{1 \leq i \leq \theta}, (\chi_i)_{1 \leq i
\leq \theta})$ be a YD-datum. Let $X$ be a vector space with basis
$x_1,\dots, x_{\theta}$. Then $\D$ defines {\rt on $X$} a Yetter-Drinfeld
module structure {\rt over $\ku\G$}  where for all $i \in \I$, $g \in \G,$
\begin{align}
\delta(x_i) &= g_i \o x_i,\\
gx_i &= \chi_i(g)x_i.
\end{align}
Then the braiding $c=c_{X,X}$ of $X$ is given by
$$c: X \otimes X \to X \otimes X,\quad c(x_i \otimes x_j) = q_{ij} x_j \otimes x_i , 1 \leq i,j \leq \theta.$$
We identify the tensor algebra $T(X)$ with the free algebra
$\ku\langle x_1,\dots,x_{\theta}\rangle$.

\medbreak Let $\mathbb{Z}^{\I}$ be a free abelian group of rank
$\theta$ with fixed basis $\alpha_1,\dots,\alpha_{\theta}$, and
$\mathbb{N}^{\I}=\{\alpha = \sum_{i=1}^{\theta} n_i \alpha_i \mid
n_1,\dots,n_{\theta} \in \mathbb{N}\}$. The homogeneous components
of a $\mathbb{Z}^{\I}$-graded vector space $Z$ will be denoted by
$Z_{\alpha}$, $\alpha \in \mathbb{Z}^{\I}$. The tensor algebra
$T(X)$ is an $\N^{\I}$-graded algebra where each $x_i$ has degree
$\alpha_i$.

\begin{Lem}\label{LemmaNichols}
\begin{enumerate}
\item $I(X)$ is an $\N^{\I}$-graded ideal of $T(X)$,
and $\NA(X)$ is an $\N^{\I}$-graded algebra and coalgebra.\label{N5}

\item\label{N7} Let $i,j \in \I$, $i \neq j$, and assume
$q_{ij}q_{ji} = q_{ii}^{a_{ij}}$ for some integer
$a_{ij}$ with $ 0 \leq -a_{ij} < \ord(q_{ii})$
(where $1 \leq \ord(q_{ii}) \leq \infty$).
Then
\end{enumerate}
\begin{equation}\label{eqn:serre}
\ad_c(x_i)^{1 - a_{ij}}(x_j) = 0 \text{ in }\NA(X).
\end{equation}
\end{Lem}
\pf (2) is shown in \cite[3.7]{AS-cambr}. For (1) see for example
\cite[Remark  2.8]{AHS}. \epf

We extend the notion of linking parameter given in \cite{AS-crelle} for data of finite Cartan type, to the general case treated here.
\begin{Def}\label{def:linking}
Vertices $i,j \in \mathbb{I}$ are called {\em linkable} if
\begin{align}
i& \not\sim j,\label{l1}\\
g_ig_j &\neq 1,\label{l2}\\
\chi_i \chi_j &=1 \text{ (the trivial character)}.\label{l3}
\end{align}
A family $\lambda = (\lambda_{ij})_{i,j \in \mathbb{I}, {\rt i \not\sim j}}$ of elements in $\ku$ is called
a {\em linking parameter} for $\D$ if for all $i,j \in \mathbb{I}$, $i \not\sim j $,
\begin{align}
\lambda_{ij}&=0 \quad \text{ if $i,j$ are not linkable},\label{l4}\\
\lambda_{ij} &= -q_{ij} \lambda_{ji}.\label{l5}
\end{align}
Given a linking parameter $\lambda$ for $\D$, vertices $i,j \in \mathbb{I}, i \not\sim j,$ are called {\em linked} if $\lambda_{ij}\neq 0$.
\end{Def}

The next Lemma generalizes \cite[Lemma 5.6]{AS2}. We include the
proof for completeness.
\begin{Lem}\label{Lemmalinking}
\begin{enumerate}
\item Let $i,j,k,l\in \mathbb{I}$.
\begin{enumerate}
\item If $i,k$ are linkable, then $q_{ji}q_{jk}=1$, and $q_{ii} = q_{kk}^{-1}=q_{ki}=q_{ik}^{-1}$.
\item If $i,k$ and $j,l$ are linkable, then $q_{ij}q_{ji}q_{kl}q_{lk}=1$.
\end{enumerate}
\item Assume that
\begin{align}\label{condition}
q_{ij}q_{ji}\neq q_{ii}^2 \text{ for all }i,j \in \mathbb{I},i \neq j.
\end{align}
Then any vertex $i \in \mathbb{I}$ is linkable to at most one $k \in \mathbb{I}$.
\end{enumerate}
\end{Lem}
\pf (1) (a) follows easily from \eqref{l3} since $q_{ik}q_{ki}=1$
by  \eqref{l1}.  (b) Since $i,k$ and $j,l$ are linkable,
$q_{ij}=q_{il}^{-1},q_{ji}=q_{jk}^{-1},q_{kl}=q_{kj}^{-1},q_{lk}=q_{li}^{-1}$
by {\rt \eqref{l3}}. Hence
\begin{equation}\label{product}
q_{ij}q_{ji}q_{kl}q_{lk} = (q_{il}q_{li})^{-1} (q_{jk}q_{kj})^{-1}.
\end{equation}
If $i \sim l$ or $j \sim k$, then $i \not\sim j$ and $k \not\sim
l$ since by assumption $i \not\sim k$ and $j \not\sim l$. Thus the
LHS of \eqref{product} is equal to 1 by \eqref{sim}.
And if $i \not\sim l$ and $j \not\sim k$, then the
RHS of \eqref{product} is equal to 1 by \eqref{sim}. This proves the claim.\\
(2) If $i \in \mathbb{I}$ is linkable to $k \in \mathbb{I}$ and to
$l \in \mathbb{I}$, then $q_{ii}^2 q_{kl}q_{lk}=1$ by (1)(b), and
$q_{ii} = q_{kk}^{-1}$ by (1)(a). Hence $q_{kl}q_{lk}=q_{kk}^2$,
and $k=l$ by assumption. \epf

For any subset $J \subset \mathbb{I}$, let $X_J=\oplus_{j \in J}
\ku x_j \in {\YDG}$. Recall the ideal $I(V)$ in Definition \ref{DefNichols}.

Let $\lambda$ be a linking parameter for the YD-datum $\D$.

\begin{Def}\label{linkingHopf}
 \begin{align}
\Uc(\D,\lambda)&=(T(X) \#\ \ku\G)/I,
\end{align}
where $I \subset T(X)\# \ku\G$ is the ideal generated by
\begin{align}
&I(X_J)\text{ for all } J \in \mathcal{X},\\
&x_ix_j - q_{ij} x_jx_i - \lambda_{ij} (1 - g_ig_j)\text{ for all } i,j \in \mathbb{I}, i \not\sim j.
\end{align}
\end{Def}
Then $\Uc(\D,\lambda)$ is a Hopf algebra in $\YDG$ with
comultiplication given by
\begin{align}
\Delta(x_i) &=g_i \otimes x_i + x_i \otimes 1, \;1 \leq i \leq \theta,\\
\Delta(g) &= g \otimes g,\; g \in \G.
\end{align}
By abuse of language we identify $x \in X$ and $g \in \G$ with
their images in $\Uc(\D,\lambda)$.

\begin{Rem}\label{Nicholscomponents}
The ideal $I(X) \subset T(X)$ is generated by
\begin{align}
&I(X_J)\text{ for all } J \in \mathcal{X},\\
&x_ix_j - q_{ij} x_jx_i \text{ for all } i,j \in \mathbb{I}, i \not\sim j.
\end{align}
Hence $\Uc(\D,0) \cong \NA(X) \# \ku\G.$
\end{Rem}

Generalizing part of \cite[Theorem 4.3]{AS-crelle}, Masuoka proved the
following result.

\begin{Theorem}\label{Masuoka} \cite[5.2]{Ma}
Let $\mathcal{X} = \{J_1,\dots,J_{t}\}$, $t \geq 1$, $J_i \neq
J_j$ for  $i\neq j \in \mathbb{I}$. For $1 \leq l \leq t$ let
$X_l=X_{J_l}$, and let $\rho_l : \NA(X_l) \to \Uc(\D,\lambda)$ be
the canonical map induced by the inclusion $T(X_l) \subset T(X)$.
Then the linear map
\begin{align*}
\NA(X_1) \o \dots \o \NA(X_t) \o \ku\G &\to \Uc(\D,\lambda),\\
r_1 \o \dots r_t \o g &\mapsto \rho_1(r_1)\cdots\rho_t(r_t)g
\end{align*}
is bijective. \qed
\end{Theorem}
Masuoka even showed that the isomorphism in Theorem \ref{Masuoka}
is a coalgebra isomorphism inducing a Hopf algebra structure on
$$(\NA(X_1) \o \dots \o \NA(X_t))  \#\ku\G$$ which is a
2-cocycle deformation. In fact he only assumes that the $\NA(X_l)$
are pre-Nichols algebras -- satisfying \eqref{N1}, \eqref{N2},
in Lemma \ref{LemmaNichols0} -- and  uses a more general equivalence relation.

\smallbreak
We close this section with a technical lemma that will be used
later.

\begin{Lem}\label{technical}
Let $j,i_1,\dots,i_n \in \I, n \geq 1, j \not\sim i_1,\dots,j
\not\sim i_n.$ Then in $\Uc(\D,\lambda)$,
\begin{align*}
x_jx_{i_1}\cdots x_{i_{n}}&= q_{i_1j}\cdots q_{i_nj} x_{i_1} \cdots x_{i_n}x_j +\\
&\sum_{1\le \nu \le n}
q_{i_1j}\cdots q_{i_{\nu -1}j}x_{i_1}\cdots x_{i_{\nu-1}}\lambda_{ji_{_\nu}}(1
- g_jg_{i_{_\nu}})x_{i_{\nu+1}} \cdots x_{i_n}.
\end{align*}
\end{Lem}
\pf If $j \in \I,x \in \Uc(\D,\lambda)$, then $x_jx = \ad
g_j(x)x_j + \ad x_j(x)$ by the definition of $\ad$; indeed,
$\Delta(x_j) = g_j \o x_j + x_j \o 1$, hence $\Ss(x_j)= -
g_j^{-1}x_j$. We apply this formula with $x = x_{i_1} \cdots
x_{i_{n}}$ and obtain
\begin{align*}
x_{j}x_{i_1} \cdots x_{i_{n}} &=
q_{i_1j} \cdots q_{i_nj}x_{i_1} \cdots x_{i_n} +\ad x_{j} (x_{i_1}\dots x_{i_n}).
\end{align*}
Now
\begin{align*}
\ad x_{j} (x_{i_1}\cdots & x_{i_n}) = \ad x_{j}\sw1(x_{i_1}) \cdots \ad x_{j(n)}(x_{i_n})\\
&=\sum_{1\le \nu \le n} \ad
g_{j}(x_{i_1})\cdots \ad g_{j}(x_{i_{\nu-1}}) \ad
x_{j}(x_{i_{\nu}})x_{i_{\nu+1}} \cdots x_{i_n}\\
&= \sum_{1\le \nu \le n}
{\rt q_{i_1j} \cdots q_{i_{\nu -1}j}} x_{i_1}\cdots x_{i_{\nu-1}}\lambda_{ji_{_\nu}}(1
- g_jg_{i_{_\nu}})x_{i_{\nu+1}} \cdots x_{i_n},
\end{align*}
where we used $j \not\sim i_{\nu}$ in the third equality, and the
claim follows. \epf

\section{A Levi-type theorem for pointed Hopf algebras}\label{section:Levi}

Let $\D$ be a YD-datum with linking parameter $\lambda$. We study
the situation when unlinked vertices are omitted. Let
\begin{itemize}
    \item []$\is = \{h\in \I: h \text{ is not linked}\}$;
    \item []$L$ a subset of $\is$;
    \item []$\I' = \I- L$;
    \item []$X'=X_{\I'}$;
    \item []$\D' = \D(\Gamma, (g_i)_{i \in \I'},
(\chi_i)_{i \in \I'})$;
    \item []$\approx$, the equivalence relation on $\I'$ defined by
the YD-datum $\D'$;
    \item []$ \lambda_{ij}'=\begin{cases} \lambda_{ij}, & \text{ if } i \not\sim j \\
0, & \text{ if } i \sim j,
\end{cases}
\text{ for all }i,j \in \I', i \not\approx j$.
\end{itemize}

Then $\lambda'$ is a linking parameter for $\D'$, since $\lambda$
is a linking parameter for $\D$. The inclusion $\iota: X'\to X$
has a section $\pi:X \to X'$ in {\rt $\YDG$} defined by
$$x_i\overset{\pi}\longmapsto x_i, \quad, x_h\overset{\pi}\longmapsto
0, \quad i\in \I',\,h\in\, L.$$

Our next Theorem can be viewed as a ``quantum version'' of the
classical Levi theorem for Lie algebras, see for instance \cite[1.6.9]{D}. We shall investigate the
case when $L = \is$ in the next section.

\begin{Theorem}\label{Levi}
The maps $\pi$ and $\iota$ induce Hopf algebra morphisms\newline
$\Phi :\Uc(\D,\lambda) \to \Uc(\D',\lambda')$ and $\Psi :
\Uc(\D',\lambda') \to \Uc(\D,\lambda)$ with $\Phi \Psi = \id$.
Then $K = \Uc(\D,\lambda)^{\co \Phi}$ is a braided Hopf algebra in
$\ydup$ and there is an isomorphism
$$K \# \Uc(\D',\lambda') \cong \Uc(\D,\lambda),$$
given by multiplication. Furthermore, the algebra $K$ is generated
by the set
\begin{equation*}\label{brackets}
S=\{\ad (x_{i_1}\cdots x_{i_n})(x_{h}) \mid h \in L, \,n \geq 0
,\, i_{\nu} \in \I',\,  i_{\nu} \sim h,0\le \nu \le n\}.
\end{equation*}
\end{Theorem}

\pf Let $\Uc = \Uc(\D,\lambda)$, $\Uc' = \Uc(\D',\lambda')$ for
brevity.

\smallbreak \emph{Existence of ${\rt \Psi}$.} We have to show that the
inclusion $\iota : T(X') \# \ku\G \hookrightarrow T(X) \# \ku\G$
maps the relations of $\Uc'$ to the relations of $\Uc$. Let $J'$
be an equivalence class of $\approx$. Then there is exactly one
equivalence class $J$ of $\sim$ with $J' \subset J$. By Lemma
\ref{LemmaNichols0}, $\iota(I(X_{J'})) \subset I(X_J)$. Let $i,j
\in \I', i \not \approx j$. We have to show that $\iota(x_ix_j
-q_{ij} - \lambda_{ij}' (1-g_ig_j))$ is a relation of $\Uc$. This
is clear if $i \not\sim j$ since $\lambda_{ij}'=\lambda_{ij}$ in
this case. If $i \sim j$ then $\lambda_{ij}'=0$ by definition, and
the relation $x_ix_j - q_{ij}x_jx_i=0$ holds in $\Uc$ by
\eqref{eqn:serre} since $q_{ij}q_{ji}=1$ follows from $i
\not\approx j$.

\smallbreak \emph{Existence of ${\rt \Phi}$.} Now we show that the
projection $\pi : T(X) \# \ku\G \to T(X') \#\ku\G$ preserves the
relations. Let $J$ be an equivalence class of $\sim$ in $\I$, and
$f \in I(X_J) \subset T(X)$; we may assume that $f$ is
$\N^{\I}$-homogeneous by Lemma \ref{LemmaNichols} \eqref{N5}. If
$f$ does not contain any variable $x_h$, $h \in L$, then
$\pi(f)=f$. Hence $f \in I(X')$ by Lemma \ref{LemmaNichols0}. Thus
$f \in T(X_J)$ is contained in the ideal $I(X')$ of $T(X')$; but
this is generated by elements in $I(X_{J'})$, $J'$ an equivalence
class of $\approx$, and elements of the form $x_ix_j -
q_{ij}x_jx_i$, where $i,j \in \I', i \not\approx j,$ and where we
can assume that $i,j \in J$. Since $\lambda_{ij}'=0$ for all $i,j
\in \I', i \not\approx j, i \sim j,$ it follows that $f$ is in the
ideal generated by the defining relations of $\Uc'$. If $f$ does
contain a variable $x_h$, where $h \in L$, then $\pi(f)=0$ since
$\pi(x_h)=0$.

\smallbreak Finally, let $i,j \in \I, i \not\sim j$. If $i \in L$
or $j \in L$, then $$\pi(x_ix_j-q_{ij}x_jx_i -
\lambda_{ij}(1-g_ig_j)) = 0,$$ since $\lambda_{ij}=0$ (because no
vertex in $L$ is linked). If $i \not\in L,j \not\in L$, then $i
\not\approx j$, and $\lambda_{ij}=\lambda_{ij}'$. Hence
$\pi(x_ix_j-q_{ij}x_jx_i -
\lambda_{ij}(1-g_ig_j))=x_ix_j-q_{ij}x_jx_i -
\lambda'_{ij}(1-g_ig_j)$ is a relation of $\Uc'$.

\smallbreak Since $\Phi \Psi = \id$ (because this holds on the generators), the multiplication map $\mu :
K \# \Uc' \to \Uc$ is an isomorphism. Let $\widetilde{K}$ be the
subalgebra of $K$ generated by $S$. Suppose we have shown that
$\widetilde{K} \Uc'$ is a subalgebra of $\Uc$. Then $\widetilde{K}
\Uc' = \Uc$ since $\widetilde{K} \Uc'$ contains the generators $g
\in \G$ and $x_i,i \in \I$ of the algebra $\Uc$. Since $\mu$ is
bijective, $\widetilde{K} = K$.

\smallbreak To prove that $\widetilde{K} \Uc'$ is a subalgebra of
$\Uc$, we have to show that
\begin{align}\label{p1}
x_j \widetilde{K} \subset \widetilde{K} \Uc' \text{ for all } j
\in \I'.
\end{align}
Then the claim follows easily by induction since the elements
$x_{j_1}\dots x_{j_n}g$, $j_1,\dots,j_n \in \I', n \geq 0, g \in
\G$, generate $\Uc '$ as a vector space, and $g \widetilde{K} =
\widetilde{K}g$. To prove \eqref{p1} it is enough to show that
\begin{align}\label{p2}
x_j \ad(x_{i_1} \cdots x_{i_n})(x_h) \in \widetilde {K} \Uc'
\end{align}
for all  $j \in \I'$, $i_1,\dots,i_n \in \I$, $n \geq 0$, $h \in
L$, $i_1 {\rt \sim} h, \dots,i_n {\rt \sim} h$. Let $x = \ad(x_{i_1}
\cdots x_{i_n})(x_h)$. {\rt By the beginning of the proof of Lemma \ref{technical}},
$$x_jx= q_{i_1j} \cdots q_{i_nj}q_{h j} xx_j + \ad x_j(x),$$
and it remains to show that $\ad x_j(x)= \ad(x_jx_{i_1} \cdots
x_{i_n})(x_h) \in \widetilde{K}\Uc'$. This is clear by definition
of $S$ if $j {\rt \sim} h$. If $j {\rt \not\sim} h$, then $\lambda_{jh}=0$
since $h \in L$, and $\ad x_j(x_h)=0$. By Lemma \ref{technical},
\begin{align*}
\ad (x_j&x_{i_1} \cdots x_{i_n})(x_h)=\ad(q_{i_1j}\cdots q_{i_nj} x_{i_1} \cdots x_{i_n}) \ad x_j(x_h)\\
 &+\ad\Big(\sum_{1\le \nu \le n}
q_{i_1j}\cdots q_{i_{\nu -1}j}x_{i_1}\cdots
x_{i_{\nu-1}}\lambda_{ji_{_\nu}}(1
- g_jg_{i_{_\nu}})x_{i_{\nu+1}} \dots x_{i_n}\Big)(x_h)\\
&=\ad\Big(\sum_{1\le \nu \le n} q_{i_1j}\cdots q_{i_{\nu
-1}j}x_{i_1}\cdots x_{i_{\nu-1}}\lambda_{ji_{_\nu}}(1 -
g_jg_{i_{_\nu}})x_{i_{\nu+1}} \dots x_{i_n}\Big)(x_h)
\end{align*}
is a $\ku$-linear combination of elements in $S$. \epf

\begin{Cor}\label{M}
In the situation of Theorem \ref{Levi}, let $\Uc=\Uc(\D,\lambda)$,
$\Uc'=\Uc(\D',\lambda')$  and $$M = \Uc/(\Uc\Uc'^+ +
\Uc(K^+)^2),$$ where $\Uc^+$ and $K^+$ are the augmentation ideals
with respect to the counit $\varepsilon$.
 Then  $x_h M \neq 0$ for any $h \in L$.
\end{Cor}
\pf By Theorem \ref{Levi}, the multiplication map $K \otimes \Uc' \to
\Uc$ is bijective; let $\psi : \Uc \to K \otimes \Uc'$ be its
inverse and
$$\varphi : \Uc \xrightarrow{\psi} K \otimes \Uc' \xrightarrow{\id \otimes \varepsilon} K.$$
Note that $\Uc K^+=K^+\Uc$, since $K=\Uc^{\co \Phi}$ and the
antipode $\Ss$ of $\Uc$ is bijective. For, $K^+$ is a submodule under  $\ad_l$ and $\ad'_r$, where
$\ad'_r(u) (x) = \Ss^{-1}(u\sw2)xu\sw1$, $x,u\in \Uc$; and the formulas
$$
ux = (\ad_l(u\_1) (x)) u\_2, \qquad xu = u\_2(\ad'_r(u\_1) (x)),
$$
hold for $x,u\in \Uc$.
Assume $x_{h}M = 0$. Then $x_{h} \in \Uc\Uc'^+ + \Uc(K^+)^2 =
K\Uc'^+ +(K^+)^2\Uc'$. Since $x_h \in K$,  it follows that
$$x_{h}=\varphi(x_{h}) \in (K^+)^2.$$
Thus by Theorem \ref{Levi}, $x_{h}$ is the $\ku$-span of products with
at least two factors of the form $\ad x_{i_1}\cdots\ad
x_{i_n}(x_{h})$, $n \geq 0$, $i_1,\dots,i_n \in J$, where $J$ is
the connected component containing $h$. Since the Nichols algebra
$\NA(V_J)$ of $V_J=\oplus_{i \in J} \ku x_i$ can be identified
with the subalgebra of $U$ generated by the elements $x_i, i \in
J$ by Theorem \ref{Masuoka}, the element $x_{h} \in \NA(V_J)$ has
degree $\geq 2$ which is impossible. \epf

\section{Perfect linkings and reduced data}\label{sect:perfect}

The goal of this section  is to study a class of pointed Hopf algebras that resembles the quantized enveloping algebras $U_q(\g)$.

\begin{Def} A linking parameter of a YD-datum $\D$ is \emph{perfect} if and only if any
vertex is linked.
\end{Def}

By Theorem \ref{Levi} for any linking parameter $\lambda$ the Hopf
algebra $\Uc(\D,\lambda)$ has a natural quotient Hopf algebra
$\Uc(\D',\lambda')$ with perfect linking parameter $\lambda'$.
This is the special case where $L = \I^s$ is the set of all
non-linked vertices.

\subsection{Reduced data}

We begin {\rt with } an alternative presentation of the Hopf algebra $\Uc(\D,\lambda)$ with perfect linking parameter;
this stresses the similarity with quantized enveloping algebras.

\begin{Def}\label{DefreducedYD-datum}
A {\em reduced YD-datum}
$$\D_{red} =\D(\G,(L_i)_{1 \leq i \leq \theta}, (K_i)_{1 \leq i \leq \theta}, (\chi_i)_{1 \leq i \leq \theta})$$
consists of an abelian group $\G$, a positive integer $\theta$,  and elements $K_i,L_i \in \G, \chi_i \in \widehat{\G}$ for all $ 1 \leq i \leq \theta$ satisfying
\begin{align}
\chi_j(K_i)&=\chi_i(L_j) \text{ for all } 1 \leq i,j \leq \theta,\label{reduced1}\\
K_iL_i &\neq 1 \text{ for all }1 \leq i \leq \theta.\label{reduced2}
\end{align}
A reduced YD-datum $\D_{red}$ is called {\em generic} if for all
$1 \leq i \leq \theta$, $\chi_i(K_i)$ is not a root of unity.
A {\em linking parameter} $\ell$ for a reduced YD-datum $\D_{red}$
is a family $\ell = (\ell_i)_{1 \leq i \leq \theta}$ of non-zero
elements in $\ku$.
\end{Def}

\begin{Def}\label{Unondegenerate}
Let $\D_{red} =\D(\G,(L_i)_{1 \leq i \leq \theta}, (K_i)_{1 \leq i
\leq \theta}, (\chi_i)_{1 \leq i \leq \theta})$ be a reduced
YD-datum with linking parameter $\ell=(\ell_i)_{1 \leq i \leq
\theta}$. Let
\begin{align}
V &= \oplus_{i=1}^{\theta}\ku v_i \in {\YDG} \text{ with basis }v_i \in V_{K_i}^{\chi_i}, 1 \leq i \leq \theta,\\
W &= \oplus_{i=1}^{\theta} \ku w_i \in {\YDG} \text{ with basis }w_i \in W_{L_i}^{\chi_i^{-1}}, 1 \leq i \leq \theta.
\end{align}
Then we define $\Uc(\D_{red}, \ell)$ as the quotient Hopf algebra
of the biproduct $T(V \oplus W) \# \ku\G$ modulo the ideal
generated by
\begin{align}
&I(V),\label{U1}\\
&I(W),\label{U2}\\
&v_iw_{j} - \chi_j^{-1}(K_i) w_{j}v_i - \delta_{ij}\ell_i(K_iL_i -
1) \text{ for all } 1 \leq i,j \leq \theta.\label{reducedlinking}
\end{align}

To a reduced YD-datum $\D_{red}$ with linking parameter $\ell$ we
associate a YD-datum $\widetilde{\D_{red}}$ and a linking
parameter $\widetilde{\ell}$ for $\widetilde{\D_{red}}$ by
\begin{align}\label{associate}
\widetilde{\D_{red}}&=\D(\G, (\widetilde{g_i})_{1 \leq i \leq 2\theta},(\widetilde{\chi_i})_{1 \leq i \leq 2\theta}), \text{ where }\\
(\widetilde{g}_1, \dots, \widetilde{g}_{2\theta}) &= (L_1,\dots,L_{\theta},K_1,\dots,K_{\theta}),\\
(\widetilde{\chi}_1, \dots, \widetilde{\chi}_{2\theta}) &= (\chi_1^{-1},\dots,\chi_{\theta}^{-1},\chi_1,\dots,\chi_{\theta}),\\
\widetilde{\ell}_{i+ \theta \,j} &= - \delta_{ij}\ell_i \text{ for
all }1 \leq i,j \leq \theta,\\ \widetilde{\ell}_{kl} &= 0 \text{
for all }1 \leq k,l  \leq 2 \theta, k \not\approx l, k >l.
\end{align}
Here $\approx$ denotes the equivalence relation of
$\widetilde{\D_{red}}.$ Note that by \eqref{l5} it suffices to
define a linking parameter $(\widetilde{\ell}_{kl})$ for all
$k>l$. Let $\widetilde{q}_{kl} =
\widetilde{\chi}_l(\widetilde{g}_k)$ for all $1 \leq k,l \leq
2\theta$, and $q_{ij}=\chi_j(K_i)$ for all $ 1 \leq i,j \leq
\theta$\label{qij}. Then  it follows from \eqref{reduced1} that for all $1
\leq i,j \leq \theta$,
\begin{align*}
\widetilde{q}_{ij}\widetilde{q}_{ji}&=(q_{ij}q_{ji})^{-1},\\
\widetilde{q}_{\theta+i, \theta+j} \widetilde{q}_{\theta+j, \theta+i}&=q_{ij}q_{ji},\\
\widetilde{q}_{i, \theta+j} \widetilde{q}_{\theta+j, i}&=1.
\end{align*}
In particular, $i \not\approx \theta+j$ for all $1 \leq i,j \leq
\theta$. Since $K_iL_i \neq1$ by  \eqref{reduced2}, it follows
that $\widetilde{\ell}$ is a linking parameter for
$\widetilde{\D_{red}}$ .

\begin{Lem}\label{Lemmareduced}
Let $\D_{red} =\D(\G,(L_i)_{1 \leq i \leq \theta}, (K_i)_{1 \leq i
\leq \theta}, (\chi_i)_{1 \leq i \leq \theta})$ be a reduced
YD-datum with linking parameter $\ell$. Then $$\Uc(\D_{red},\ell)
\cong \Uc(\widetilde{\D_{red}}, \widetilde{\ell}).$$
\end{Lem}
\pf This follows from the defining relations using Remark
\ref{Nicholscomponents}. \epf

\begin{Lem}\label{perfectreduced}
Let $\D=\D(\G,(g_i)_{1 \leq i \leq \theta}, (\chi_i)_{1 \leq i
\leq \theta})$ be a YD-datum satisfying \eqref{condition}, and let
$\lambda$ be a perfect linking parameter for $\D$. Then there {\rt is}
a reduced YD-datum $\D_{red}$ and a linking parameter $\ell$ for $
\D_{red}$ such that
$$\Uc(\D,\lambda) \cong \Uc(\D_{red},\ell)$$
as Hopf algebras.
\end{Lem}
\pf If $i\in \I$ then by Lemma \ref{Lemmalinking} (2) there
exists a unique $i^0\in \I$ such that $i$ and $i^0$ are linked.
Thus $\I \to \I,\;i\mapsto i^0,$ is an involution on the set of
vertices. By Lemma \ref{Lemmalinking} (1)(b),
$q_{ij}q_{ji}q_{i^0j^0}q_{j^0i^0} =1$ for all $i,j \in \I.$ Hence
$$\X \to \X,\quad J \mapsto J^0=\{j^0 \mid j \in J\},$$
is an involution on the set of equivalence classes, and $J \cap
J^0 =\emptyset$ for all $J \in \X$ since $i \not\sim i^0$ for all
$i \in \I$. Therefore after renumbering the indices we may assume
that $\I = \I^- \cup \I^+,$ where $\I^-=\{1,\dots,\theta_1\}$,
$\I^+=\{\theta_1+1,\dots,2\theta_1\}$, $\theta=2\theta_1$, and
$i^0=\theta_1+i$ for all $i \in \I^-$. Moreover there are subsets
$\X^-,\X^+ \subset \X, \X^+=\{ J^0 \mid J \in \X^-\}$ such that

$$\I^- = \bigcup_{J \in  \X^-} J, \quad \I^+ = \bigcup_{J \in  \X^+} J.$$
Then for all $1 \leq i \leq \theta_1$, $\chi_i = \chi_{\theta_1
+i}^{-1}$, and $\widetilde{\ell}_{\theta_1 +i,j} \neq 0$ since $i,
\theta_1 +i$ are linked. Define $\D_{red}(\Gamma, (L_i)_{1 \leq
\theta_1},(K_i)_{1 \leq i \leq \theta_1}, (\eta_i)_{1 \leq i \leq
\theta})$ and $\ell=(\ell_i)_{1 \leq i \leq \theta_1}$  by
\begin{align*}
(g_1,\dots,g_{2\theta_1}) &= (L_1,\dots,L_{\theta_1},K_1,\dots,K_{\theta_1}),\\
(\chi_1,\dots,\chi_{2\theta_1})&=(\eta_1^{-1},\dots,\eta_{\theta_1}^{-1},\eta_1,\dots,\eta_{\theta_1}),\\
\ell_i &= -\lambda_{\theta_1 +i,i}, \qquad  1 \leq i \leq
\theta_1.
\end{align*} Then the lemma follows from Lemma \ref{Lemmareduced}
since $\D=\widetilde{\D_{red}}$, $\widetilde{\ell}= \lambda$. \epf

\subsection{The Hopf algebra $\Uc(\D_{red},\ell)$ as a quotient of a Drinfeld double}

 \

For the rest of this section we fix a reduced YD-datum
$$\D_{red} =\D(\G,(L_i)_{1 \leq i \leq \theta}, (K_i)_{1 \leq i \leq \theta}, (\chi_i)_{1 \leq i \leq \theta})$$ with
linking parameter $\ell = (\ell_i)_{1 \leq i \leq \theta}$, and
denote
$\U=\Uc(\D_{red},\ell)$. We shall describe $\U$ as a quotient of a Drinfeld double.

\smallbreak
The images of $v_i$ and $w_i$ in $\U$ will again be denoted by $v_i$ and $w_i$. Let
\begin{align}
E_i=v_{i}, \quad F_i=w_iL_{i}^{-1} \text{ in } \U, \quad 1 \leq i \leq \theta,
\end{align}
and let $\U^-$ (resp. $\U^+$) be the subalgebra of $\U$ generated by $F_1,\dots,F_{\theta}$ (resp. $E_1,\dots,E_{\theta}$). Then
\begin{align}
\label{eqn:g-con-ei}gE_ig^{-1}&=\chi_i(g)E_i,\\
\label{eqn:g-con-fi}gF_ig^{-1} &= \chi_i^{-1}(g)F_i,\\
E_iF_i - F_i E_i &= \delta_{ij}\ell_i(K_i -L_i^{-1}),\\
\label{eqn:comul-ei}\Delta(E_i) &= K_i \otimes E_i + E_i \otimes 1,\\
\label{eqn:comul-fi}\Delta(F_i) &=1 \otimes F_i + F_i \otimes
L_i^{-1}
\end{align}
in $\U$, for all $1 \leq i \leq \theta$, $g\in \Gamma$.
\end{Def}

\begin{Rem}\label{relationsF}
The Hopf algebra $\Uc(\D_{red},\ell)$ does not depend on the
choice of the non-zero scalars $\ell_i$. By rescaling the
variables $v_i$ we could assume that $\ell_i=1$ for all $i$. By
the same reason we could assume that the only values
$\lambda_{ij}$ of a linking parameter for a YD-datum $\D$ are 0 or
1.

\medbreak Since $\Ss(F_i)=-F_iL_i=-w_i$ in $\NA(W) \# \ku\G$ for
all $1 \leq i \leq \theta$, the relations of the elements $w_i$ in
$\Uc(\D_{red},\ell)$ may be equivalently expressed by the
following relations in the $F_i$. If $f$ is an element of the free
algebra in the variables $x_1,\dots,x_{\theta}$, and the relation
$f(w_1,\dots,w_{\theta})=0$ holds in $\Uc(\D_{red},\ell)$, then
$\widetilde{f}(F_1,\dots,F_{\theta})=0$, where
$$\ku\langle x_1,\dots,x_{\theta}\rangle \to \ku\langle x_1,\dots,x_{\theta}\rangle,
\quad f \mapsto \widetilde{f},$$ is the vector space isomorphism
mapping a monomial $x_{i_1}x_{i_2} \cdots x_{i_n}$ onto $(-1)^n
x_{i_n} \cdots x_{i_2}x_{i_1}$.
\end{Rem}


\smallbreak Let $\Lambda$ be the free abelian group with basis
$z_1,\dots,z_\theta,$ and define characters $\eta_j \in
\widehat{\Lambda}$ by $\eta_j(z_i) = \chi_j^{-1}(L_i)$, $1 \leq i,j
\leq \theta$. Then $W$ is Yetter-Drinfeld module in $\YDL$ with
$w_i \in W_{z_i}^{\eta_i}, 1 \leq i \leq \theta,$ and $W$ has the
same braiding as an object in ${\YDG}$ or in ${\YDL}$. Hence
$\NA(W)$ is a Hopf algebra in ${\YDG}$ and ${\YDL}$. Let
$A=\mathfrak{B}(V) \# k[\G]$ and $U= \mathfrak{B}(W) \#
k[\Lambda]$.

\smallbreak Generally for Hopf algebras $A$ and $U$ a linear map $\tau : U \otimes A \to \ku$ is a {\em
skew-pairing} \cite[Definition 1.3]{DT} if
\begin{align}\label{form1}
 \tau(u,aa') &= \tau(u\sw2,a) \tau(u\sw1,a'),\\
\tau(uu',a) &= \tau(u,a\sw1) \tau(u',a\sw2),\\
\tau(1,a) &= \varepsilon(a),\;\tau(u,1) = \varepsilon(u),
\end{align}
for all $u,u' \in U$ and $a,a' \in A$.

A skew-pairing $\tau$ defines a 2-cocycle $\sigma : (U \otimes A)
\otimes (U \otimes A) \to k$ by
\begin{equation}\label{tausigma}
\sigma(u \otimes a,u' \otimes a')= \varepsilon(u) \tau(u',a)
\varepsilon(a')
\end{equation}
for all $u,u' \in U,a,a' \in A$. Let $(U \otimes A)_{\sigma}$ be
the 2-cocycle twist of the tensor product Hopf algebra $U \o A$.
Thus $(U \otimes A)_{\sigma}$  coincides with $U \otimes A$ as a
coalgebra with componentwise comultiplication and its algebra
structure is defined by
\begin{align}\label{Eq2CocyleMult}
(u{\otimes}a)(u'{\otimes}a') &= \sigma(h\sw1,h'\sw1)h\sw2h'\sw2 \sigma^{-1}(h\sw3,h'\sw3)\\
&=u\tau(u'_{(1)},
a_{(1)})u'_{(2)}{\otimes}a_{(2)}\tau^{-1}(u'_{(3)}, a_{(3)})a'\notag
\end{align}
for all $u, u' \in U,a, a' \in A$.
Note that
$$\tau^{-1}(u,a) = \tau(\Ss(u),a)= \tau(u,\Ss^{-1}(a))$$
for all $u \in U, a\in A$.

\medbreak Part (1) of the next result is a special case of \cite[Theorem 8.3,
Corollary 9.1]{RS1}, part (2) is shown in \cite[Theorem 4.4]{RS2} for data
of (finite) Cartan type, and in general by similar methods in
\cite[Theorem 5.3]{Ma}.
Let $A$, $U$ be the bosonizations defined above.

\begin{Theorem}\label{form}
\emph{(1)} There is a unique skew-pairing $\tau : U \otimes A \to
k$ with
\begin{align*}
\tau(z_i,g) &=\chi_i^{-1}(g),& \tau(z_i,v_j)&=0, \\
\tau(w_i,g)&=0,&\tau(w_i,v_j)&=-\delta_{ij}\ell_i
\end{align*}
for all $1 \leq i,j \leq \theta$ and $g \in \G$.\label{form2}

\emph{(2) }Let $\sigma$ be the 2-cocycle corresponding to $\tau$
by \eqref{tausigma}. Then there is an isomorphism of Hopf algebras
\begin{equation*}\label{iso}
\U \cong  (U \otimes A)_{\sigma}/(z_i \otimes L_i^{-1} - 1 \otimes 1 \mid 1 \leq i \leq \theta),
\end{equation*}
mapping $w_i,1 \leq i \leq \theta,$ and  $v_{j},1 \leq j \leq
\theta,$ respectively $g \in \Gamma$ onto the residue classes of
$w_i \otimes 1,1 \otimes v_j$, respectively $1 \otimes g$. \qed
\end{Theorem}

The following decomposition result is a special case of \cite[Theorem
5.2]{Ma}. By definition of $\U$ there are algebra maps $\rho_V :
\NA(V) \to \U,\; \rho_W : \NA(W) \to \U$ and $\rho_{\G} : \ku\G
\to \U$, given by $\rho_V(v_i) = v_i$, $\rho_W(w_i)= w_i$,
$\rho_{\G}(g) = g$, for all $1 \leq i \leq \theta, g \in \G$.
Clearly, the image of $\rho_V$ coincides with $\U^+$-- but the
image of $\rho_W$ is not $\U^-$.

\begin{Cor}\label{decomposition}
\emph{(1)} The multiplication map
$$ \NA(V) \otimes \NA(W) \otimes \ku\G \to \U,\; v \otimes w \otimes g \mapsto \rho_V(v) \rho_W(w) \rho_{\G}(g),$$
is a coalgebra isomorphism.

\emph{(2)} The multiplication map $\U^- \otimes \U^+ \otimes k[\G]
\to \U$ is an isomorphism of vector spaces.

\end{Cor}
\pf (1) follows from Theorem \ref{Masuoka} and Lemma
\ref{perfectreduced}. We prove (2). By (1) we may
identify $\NA(V),\NA(W)$ and $\ku\G$ with subalgebras of $\U$. We
first claim that the multiplication map $\ku\G \otimes \NA(W) \to
\U$ defines an isomorphism $$\ku\G \otimes \NA(W)  \cong \ku\G
\NA(W) = \NA(W) \ku\G.$$ The multiplication map defines an
isomorphism $\NA(W) \otimes \ku\G \cong \NA(W) \ku\G$ by (1).
Since $gw_i = \chi^{-1}(g) w_ig$ for all $ 1\leq i \leq \theta, g
\in \G$, $\NA(W)$ has a vector space basis $(w_b)_{b \in B}$ such
that $gw_b= \chi_b(g) w_bg$ for all $b \in B, g \in \G$, where the
$\chi_b$ are characters of $\G$. Hence also $\ku\G \otimes \NA(W)
\to \U$ is injective, and the claim follows. Then it follows from
(1) that the multiplication map
\begin{equation}\label{newiso}
\NA(V) \otimes \ku\G \otimes \NA(W) \to \U
\end{equation}
is bijective.  By (1), $\NA(V) = \U^+$, and $\NA(V) \# {\rt \ku\G}
\cong \U^+ \ku\G$ is a Hopf subalgebra of $\U$. Also, $\Ss(F_i) =
-w_i$ for all $1 \leq i \leq \theta$, and $\Ss(\U^-)=\NA(W)$. By
\eqref{newiso} the composition
$$\U^+ \otimes \ku\G \otimes \U^-\cong \U^+ \ku\G
\otimes \U^- \xrightarrow{\Ss \otimes \Ss} \U^+\ku\G\otimes \NA(W)
\cong \U,$$ mapping $x \otimes g \otimes y$ onto $\Ss(yxg)$ for
all $x \in \U^+, g \in \G, y \in \U^-$ is bijective. Thus
multiplication defines an isomorphism $\U^- \otimes \U^+ \otimes
\ku\G \to \U$. \epf

\subsection{A bilinear form}

\

We will now see that the form $\tau : U \otimes A \to \ku$ defines
in a natural way a form $(\;,\;) : \U^- \otimes \U^+ \to \ku$.
This is the form we will use later on. {\rt  It plays the same role as Lusztig's form $(\;,\;) : {\bf f} \ot {\bf f} \to \mathbb{Q}(v)$.}

\smallbreak Let $\pi_{\G} : \mathfrak{B}(V) \# k[\G] \to k[\G]$ be
the projection defined by $\pi_{\G}(x \otimes g) =
\varepsilon(x)g$ for $x \in \NA(V)$, $g \in \G$. Clearly,
$\NA(V)=A^{\co \pi_{\G}}$.  Let $\pi_{\Lambda}: \mathfrak{B}(W) \#
k[\Lambda] \to k[\Lambda]$ be the analogous  projection of $U$ to
$k[\Lambda]$.  We have $$\Ss(\NA(W)) = \Ss(U^{\co \pi_{\Lambda}}) =
{^{\co\pi_{\Lambda}}U},$$ see Subsection
\ref{subsection:bosonization}; thus ${^{\co\pi_{\Lambda}}U}$ is
generated as an algebra by the elements
$w_1z_1^{-1},\dots,w_{\theta}z_{\theta}^{-1}$ since $\Ss(w_i) =
-z_i^{-1}w_i = -q_{ii}w_i z_i^{-1}$ for all $1 \leq i \leq
\theta$.

\begin{Cor}\label{iota}
\emph{(1)} The Hopf algebra map $\varphi^+ : \NA(V) \# \ku\G \to
\U$ given by
$$ \varphi^+(v_i)
=E_i,\quad\varphi^+(g)=g, \quad 1 \leq i \leq \theta, \; g
\in\G,$$ is injective. {\rt In particular,
$$\iota^+ : \U^+ \to A^{\co \pi_{\G}},\;
\iota^+(E_i) = v_i,\quad 1 \leq i \leq \theta$$
is a well-defined algebra isomorphism.}

\emph{(2)} The Hopf algebra map $\varphi^- : \NA(W) \#
\ku\Lambda \to \U$ given by
$$\varphi^-(w_i)= F_iL_i,\quad \varphi^-(z_i) =L_i, \quad 1 \leq i \leq \theta,$$ induces a
bijection  between the subalgebras ${^{\co\pi_{\Lambda}}U}$ and
$\U^-$. In particular,
$$\iota^- : \U^- \to {^{\co\pi_{\Lambda}}U},\quad
\iota^-(F_i) = w_i z_i^{-1},1 \leq i \leq \theta,$$ is a well-defined algebra isomorphism.

{\rt \emph{(3)} Let $\kappa$ be the bijective map of Lemma \ref{leftright} with respect to the projection $\pi = \pi_{\Lambda} : \NA(W) \# \ku\Lambda \to \ku\Lambda$.
Then $\varphi^- \kappa$ defines a bijective linear map
$$\widetilde{\kappa} : \NA(W) \to \U^-,\;
\widetilde{\kappa}(w_{i_1} \cdots w_{i_n}) = \prod_{l=1}^n q_{i_li_l} \prod_{k > l} q_{i_ki_l}^{-1} F_{i_1} \cdots F_{i_{n}}$$
for all $1 \le i_1,\dots,i_n \le \theta,n \ge 1.$}
\end{Cor}

\pf (1) follows from Corollary \ref{decomposition} (1).

(2) The Hopf algebra map $\varphi^-$ is the composition of the
Hopf algebra maps $\NA(W) \# \ku\Lambda \to \NA(W) \# \ku\G$,
$w_i \mapsto w_i$, $z_i \mapsto L_i$, $1 \leq i \leq \theta$ and
$\NA(W) \# \ku\G \to \U$, $w_i \mapsto w_i = F_iL_i$, $g \mapsto
g$, $1 \leq i \leq \theta$, $g \in \G$. The restriction of
$\varphi^-$ is an isomorphism from $\mathfrak{B}(W)$ to the
subalgebra $\ku \langle w_1,\dots,w_{\theta}\rangle$ of $\U$, by (1). Hence
$\varphi^-$ induces an isomorphism between $\Ss(\mathfrak{B}(W))=
{^{\co\pi_{\Lambda}}U}$ and $\Ss(\ku\langle w_1,\dots,w_{\theta}\rangle)= \U^-$.
Its inverse is $\iota^-$.

{\rt (3) follows from  (2) and the formula for $\kappa$ in Lemma \ref{leftright}.}
\epf

\begin{Def}
The $\ku$-bilinear form $(\;,\;) : \U^- \otimes \U^+ \to \ku$ is defined by $(x,y)=\tau(\iota^-(x),\iota^+(y))$ for all $x \in \U^-,y \in \U^+$.
\end{Def}

\medskip
If $\alpha = \sum_{i=1}^{\theta} n_i \alpha_i \in
\mathbb{Z}^{\I}$, $n_1,\dots,n_{\theta} \in \mathbb{Z},$ we let
$|\alpha| = \sum_{i=1}^{\theta} n_i$, and
\begin{equation}
\chi_{\alpha} = \chi_1^{n_1} \cdots
\chi_{\theta}^{n_{\theta}},\quad K_{\alpha} = K_1^{n_1} \cdots
K_{\theta}^{n_{\theta}}, \quad L_{\alpha} = L_1^{n_1} \cdots
L_{\theta}^{n_{\theta}}.
\end{equation}

\smallbreak
The Hopf algebras $U=\mathfrak{B}(W) \# k[\Lambda]$ and
$A=\mathfrak{B}(V) \# k[\G]$ are $\mathbb{N}^{\I}$-graded as
algebras and coalgebras where the elements $v_i$  have degree
$\alpha_i$ and the elements $w_i$ have degree $-\alpha_i$ for all
$1 \leq i \leq \theta$, and the elements of the groups $\Lambda$
and $\G$ have degree 0. Hence the algebras $\U^+,\U^-$ are
$\mathbb{N}^{\I}$-graded by Corollary \ref{iota}, where  the degree of
$E_i$ is $\alpha_i$ and the degree of $F_i$ is $- \alpha_i$, for
all $1 \leq i \leq \theta$.

\medbreak We collect some important properties of the forms $\tau$
and $(\;,\;)$.

\begin{Theorem}\label{rules}
\begin{enumerate}
\item $\tau(uz,ag)=\tau(u,a)\tau(z,g)$,
for all $u \in {^{\co\pi_{\Lambda}}U}$, $a \in
A^{\co\pi_{\G}}=\mathfrak{B}(V)$, $z \in \Lambda, g \in \G$.
\label{form3}

\item For all $ \alpha, \beta \in \mathbb{N}^{\I}$,
$\alpha \neq \beta,$ $\tau(U_{-\alpha},A_{\beta})=0$.\label{form4}

\item For all $\alpha \in \mathbb{N}^{\I}$, the restriction
of $\tau$ to $\mathfrak{B}(W)_{-\alpha} \times
\mathfrak{B}(V)_{\alpha}$ is non-degenerate.\label{form5}

\item For all $ \alpha, \beta \in \mathbb{N}^{\I},\alpha \neq
\beta$, $({\U^-}_{-\alpha},{\U^+}_{\beta})=0$. \label{form6}

\item For all $\alpha \in \mathbb{N}^{\I}$, the restriction of
the form $(\;,\;)$ to $\U^-_{-\alpha} \times \U^+_{\alpha}$ is
non-degenerate.\label{form7}

\item For all $1 \leq i,j \leq \theta,
(F_i,E_j) = -\delta_{ij}\ell_i.$
\end{enumerate}
\end{Theorem}
\pf (1) The proof follows from the claims
\begin{enumerate}
\item [(a)] $\tau(z,a)=\tau(z,\pi_{\G}(a))$,
\item [(b)] $\tau(u,g) = \tau(\pi_{\Lambda}(u),g)$
\end{enumerate}
for all $z \in \Lambda$, $a \in A$, $u \in U$, $g \in \G$. For, suppose
(a) and (b) hold and $u,a,z,g$ satisfy the hypothesis of (1). Then
using Theorem \ref{form} we calculate
\begin{align*}
\tau(uz,ag)&= \tau(u,a\sw1g)\tau(z,a\sw2g)\\
&=\tau(u\sw2,a\sw1)\tau(u\sw1,g)\tau(z,a\sw2)\tau(z,g)\\
&=\tau(u\sw2,a\sw1)\tau(\pi_{\Lambda}(u\sw1),g)\tau(z,\pi_{\G}(a\sw2))\tau(z,g)\\
&= \tau(u,a) \tau(z,g).
\end{align*}
We prove (a). Since $z \in G(U)$, the map $\tau(z,-) : A \to k$ is
an algebra map by part (1) of Theorem \ref{form}. Since
$K_iv_iK_i^{-1} = q_{ii} v_i$ and $q_{ii} \neq 1$ for all $1 \leq
i \leq \theta$ it follows that any algebra map from $A$ to $k$
vanishes on each $v_i$. Thus $\tau(z,a) = \tau(z,\pi_{\G}(a))$ for
all $a \in A$. The second claim (b) follows similarly using
$\tau(-,g)$ in place of $\tau(z,-)$.

(2) follows from Theorem \ref{form} (1) and the fact that the
comultiplications of $U$ and $A$ are $\mathbb{Z}^{\I}$-graded.

(3) Since all the $\ell_i$ are non-zero, the form $\tau$ restricts
to a non-degenerate pairing between $\NA(W)$ and $\NA(V)$ (see
\cite{RS1} or \cite[Remark 3.3]{RS2}). Hence the claim in (3)
follows from (2).

(4) and (5) follow from (2) and (3) using (1), since $U= {^{\co \pi_{\Lambda}}U}\Lambda$.

(6) follows from Theorem \ref{form} (1).
\epf

\subsection{{\rt Further} properties of the bilinear form}\label{sectionquasi}

\

We now discuss some further properties of the bilinear form {\rt following \cite[Chapters 3 and 4]{L}}; in
particular, we study a universal element in some completion of
$\U$. In the case of reduced data of Cartan type, it will give
rise to Casimir elements, up to some suitable modification.

\medbreak In \cite[1.2.13]{L} Lusztig introduces two skew-derivations $r_i$ and $_ir$. We need four such maps. The comultiplication of $\U$ defines  skew-derivations
$r_i,r_i' : \U^+ \to \U^+$ and $s_i,s_i' : \U^- \to \U^-$ for all
$ 1 \leq i \leq \theta$ in the following way.

Since  $\Delta(E_i) = K_i \o E_i + E_i \o 1$, for all $1 \leq i
\leq \theta$, it follows that for all $\alpha \in \N^{\I}$ and $y
\in \U^+_{\alpha}$, $\Delta(y)$ has the form

\begin{align}
\Delta(y) &= y \o 1 + \sum_{i=1}^{\theta} r_i(y)K_i \o E_i + \text{ terms of other degrees},\label{r1}\\
\Delta(y)&= K_{\alpha} \o y + \sum_{i=1}^{\theta} E_i K_{\alpha - \alpha_i} \o r'_i(y)+ \text{ terms of other degrees},\label{r2}
\end{align}
where $r_i(y), r'_i(y)$ are uniquely determined elements in $\U^+_{\alpha - \alpha_i}$.
Degree refers to the standard $\Z^\I$-grading in the tensor product.
Then for all $y,y' \in \U^+$ and $1 \leq i \leq \theta$,
\begin{align}
r_i(yy') &= y r_i(y') + r_i(y)(K_i \on y'),\label{derivation1}\\
r'_i(yy')&= (L_i \on y) r'_i(y') + r'_i(y) y'.\label{derivation2}
\end{align}

This follows from $\Delta(yy')= \Delta(y)\Delta(y')$ by comparing
coefficients. Note that $r_i(E_j) = \delta_{ij}$,
$r'_{ij}(E_j)=\delta_{ij}$ for all $1 \leq i,j \leq \theta$.

In the same way it follows from $\Delta(F_i) = 1 \otimes F_i + F_i
\otimes L_i^{-1}$ for all $1 \leq i \leq \theta$ that for all
$\alpha \in \N^{\I}$ and $x \in \U^-_{-\alpha}$,
\begin{align}
\Delta(x) &= x \o L_{\alpha}^{-1} + \sum_{i=1}^{\theta} s_i(x) \o F_i L_{\alpha - \alpha_i}^{-1} + \text{ terms of other degrees},\label{s1}\\
\Delta(x)&= 1 \o x + \sum_{i=1}^{\theta} F_i \o s'_i(x)L_i^{-1} + \text{ terms of other degrees},\label{s2}
\end{align}
where $s_i(x)$, $s'_i(x) \in \U^-_{-\alpha + \alpha_i}$ are
uniquely determined elements. Then for all $1 \leq i,j \leq
\theta$, $s_i(F_j) = \delta_{ij}$, $s'_{ij}(F_j)=\delta_{ij}$, and
for all $x,x' \in \U^-$ and $1 \leq i \leq \theta$,
\begin{align}
s_i(xx') = (K_i^{-1}\on x) s_i(x') + s_i(x) x',\label{derivation3}\\
s'_i(xx')= x s'_i(x') + s'_i(x) (L_i^{-1} \on x').\label{derivation4}
\end{align}

\medbreak

{\rt The next Propositions \ref{rulesF}, \ref{rulesE} extend \cite[3.1.6]{L}.}

\begin{Prop}\label{rulesF}
For all $x \in \U^-,y \in \U^+$ and $1 \leq i \leq \theta$,
\begin{align*}
yF_i -F_iy &= \ell_i (r_i(y) K_i - L_i^{-1} r'_i(y)),\label{rule1}\tag{1}\\
(xF_i,y) &= (x,r_i(y)) (F_i,E_i),\label{rule2}\tag{2}\\
(F_ix,y) &= (x, r'_i(y)) (F_i,E_i).\label{rule3}\tag{3}
\end{align*}
\end{Prop}
\pf (1) The function $d_i : \U^+ \to \U,\; y \mapsto r_i(y) K_i -
L_i^{-1} r'_i(y),$ is a derivation since for all $y,y' \in \U^+$,
\begin{align*}
d_i(yy')&=r_i(yy') K_i - L_i^{-1} r'_i(yy')\\
&= y r_i(y')K_i + r_i(y)(K_i \on y')K_i - L_i^{-1}(L_i \on y) r'_i(y') - L_i^{-1} r'_i(y) y'\\
&=y r_i(y')K_i + r_i(y) K_i y' - y L_i^{-1} r'_i(y') - L_i^{-1} r'_i(y) y'\\
&= d_i(y)y' + yd_i(y'),
\end{align*}
where we have used \eqref{derivation1}, \eqref{derivation2} and the equalities
$$(K_i \on y')K_i= K_i y',\quad L_i^{-1}(L_i \on y)=yL_i^{-1}.$$
Moreover, $d_i(E_j) = \delta_{ij}(K_i - L_i^{-1})$, for all $1
\leq j \leq \theta$. Since both sides of (1) are derivations
having the same values on the generators $E_j$ of $\U^+$, the
claim follows.

(2) We can assume that $y \in \U^+_{\alpha}$, where $\alpha \in
\N^{\I}$. Let $u=\iota^-(x)$, $a =\iota^+(y)$. Then
$$\Delta(a) = a \o 1 + \sum_{i=1}^{\theta} \bar{r}_i(a)K_i \o v_i + \text{ terms of other degrees},$$
where $\bar{r}_i(a) = \iota^+(r_i(y))$ by Corollary \ref{iota}
(1). Hence, by Lemmas \ref{form}, \ref{rules},
\begin{align*}
(xF_i,y)&= \tau(uw_iz_i^{-1}, a)\\
&= \tau(u,a\sw1) \tau(w_iz_i^{-1},a\sw2)\\
&= \tau(u,\bar{r}_i(a)K_i) \tau(w_iz_i^{-1},v_i)\\
&= \tau(u,\bar{r}_i(a)) \tau(w_iz_i^{-1},v_i)\\
&=(x,r_i(y)) (F_i,E_i).
\end{align*}

(3) is proved in the same way as (2).
\epf

\begin{Prop}\label{rulesE}
For all $x \in \U^-,y \in \U^+$ and $1 \leq i \leq \theta$,
\begin{align*}
E_ix -  xE_i &= \ell_i(K_i s_i(x) - s'_i(x) L_i^{-1}),\label{rule4}\tag{1}\\
(x,E_iy) &= (s_i(x),y) (F_i,E_i),\label{rule5}\tag{2}\\
(x,yE_i) &= (s'_i(x),y) (F_i,E_i).\label{rule6}\tag{3}
\end{align*}
\end{Prop}
\pf Similar to the proof of Proposition \ref{rulesF} using Corollary
\ref{iota} (2). \epf

Recall that the form $(\phantom{a},\phantom{a}) : \U^-_{-\alpha}
\times \U^+_{\alpha} \to k$ is non-degenerate by Theorem \ref{rules}
\eqref{form7}, for all $\alpha \in \N^{\I}$.

\begin{Def} For all $\alpha \in \N^{\I}$, let $x_{\alpha}^k$,
$1 \leq k \leq d_{\alpha} =\dim \U^-_{-\alpha}$, be a basis of
$\U^-_{-\alpha}$; and $y_{\alpha}^k, 1 \leq k \leq d_{\alpha},$
the dual basis of $\U^+_{\alpha}$ with respect to
$(\phantom{a},\phantom{a})$. Define
$$\theta_{\alpha} =\sum_{k=1}^{d_{\alpha}} x_{\alpha}^k \o y_{\alpha}^k.$$
We set $\theta_{\alpha}=0$ for all $\alpha \in \Z^\I$ and $\alpha
\notin \N^{\I}$. The following  formal element is instrumental to
the definition of the quantum Casimir element:

$$\Omega= \sum_{\alpha \in \N^{\I}}\sum_{k=1}^{d_{\alpha}} \Ss(x_{\alpha}^k)y_{\alpha}^k.$$
\end{Def}

We collect some general properties of the family
$(\theta_{\alpha})$ generalizing \cite[4.2.5]{L}.

\begin{Theorem}\label{quasi}
Let $\alpha \in \N^{\I}$ and $1 \leq i \leq \theta$. Then in $\U \o
\U$,
\begin{align*}
(E_i \o 1) \theta_{\alpha} + (K_i \o E_i) \theta_{\alpha - \alpha_i} &= \theta_{\alpha}(E_i \o 1) + \theta_{\alpha - \alpha_i}(L_i^{-1} \o E_i),\tag{1}\\
(1 \o F_i)\theta_{\alpha} + (F_i \o L_i^{-1}) \theta_{\alpha - \alpha_i} &= \theta_{\alpha}(1 \o F_i) + \theta_{\alpha - \alpha_i}(F_i \o K_i).\tag{2}
\end{align*}
\end{Theorem}
\pf Both equalities hold when $\alpha - \alpha_i \not\in \N^{\I}$
since then $E_i$ commutes with the elements $x^k_{\alpha}$ which
are products of $F_j's$ where $j \neq i$, and similarly $F_i$
commutes with the elements $y^k_{\alpha}$.

By definition the equality in (1) means  that
\begin{align*}
&\sum_k E_i x_{\alpha}^k \o y_{\alpha}^k + \sum_l K_i x_{\alpha -
\alpha_i}^l \o E_i y_{\alpha - \alpha_i}^l \\ &-\sum_k
x_{\alpha}^kE_i \o y_{\alpha}^k - \sum_l x_{\alpha - \alpha_i}^l
L_i^{-1} \o y_{\alpha - \alpha_i}^l E_i=0\end{align*} in $\U
\otimes \U^+_{\alpha}$, or equivalently, by non-degeneracy of
$(\phantom{a},\phantom{a})$, that
\begin{align*}&\sum_k (E_i x_{\alpha}^k
-x_{\alpha}^kE_i)(z,y_{\alpha}^k)
+ \sum_l K_i x_{\alpha - \alpha_i}^l(z,E_i y_{\alpha - \alpha_i}^l)\\
 -&\sum_l x_{\alpha - \alpha_i}^l L_i^{-1}(z,y_{\alpha - \alpha_i}^l
 E_i)=0,
\end{align*}
for all $z \in \U_{-\alpha}^-$. Now we apply Proposition \ref{rulesE}
\eqref{rule4}, \eqref{rule5} and \eqref{rule6} to the summands of
the first, second and third sum, collect coefficients of $K_i$ and
$L_i^{-1}$ and obtain the following equivalent form of (1)
\begin{align*}
&K_i\left((F_i,E_i)\sum_l x_{\alpha - \alpha_i}^l (s_i(z),y_{\alpha - \alpha_i}^l) + \ell_i\sum_k s_i(x_{\alpha}^k)(z,y_{\alpha}^k)\right)\\
-&\left(\ell_i\sum_ks_i'(x_{\alpha}^k)(z,y_{\alpha}^k) +
(F_i,E_i)\sum_l x_{\alpha - \alpha_i}^l (s_i'(z), y_{\alpha -
\alpha_i}^l)\right)L_i^{-1} =0.
\end{align*}
Since the tensorands of $\theta_{\alpha}$ and $\theta_{\alpha - \alpha_i}$ are dual bases, we see that
\begin{align*}
&\sum_l x_{\alpha - \alpha_i}^l (s_i(z),y_{\alpha - \alpha_i}^l) =s_i(z),\\
&\sum_k s_i(x_{\alpha}^k)(z,y_{\alpha}^k))= s_i(\sum_kx_{\alpha}^k(z,y_{\alpha}^k))=s_i(z).
\end{align*}
Since $(F_i,E_i)=-\ell_i$, it follows that the coefficient of
$K_i$ is zero. Similarly the coefficient of $L_i^{-1}$ is zero
since
\begin{align*}
&\sum_ks_i'(x_{\alpha}^k)(z,y_{\alpha}^k)= s_i'(z),\\
&\sum_l x_{\alpha - \alpha_i}^l (s_i'(z), y_{\alpha - \alpha_i}^l)= s_i'(z).
\end{align*}
(2) is proved in the same way using Proposition \ref{rulesF} instead of
Proposition \ref{rulesE}. \epf

\subsection{Data of Cartan type}

\

Let $(a_{ij})_{1 \leq i,j \leq \theta}$ be a generalized Cartan
matrix, that is, $(a_{ij})_{1 \leq i,j \leq \theta}$ is a matrix
with has  integer entries such  $a_{ii} =2$ for all $1 \leq i \leq
\theta$, and  for all $1 \leq i,j \leq \theta,i \neq j$, $a_{ij}
\leq 0$, and if $a_{ij}=0$, then $a_{ji} =0$.

\begin{Def}\label{CartanYD-datum}
Let $\D = \D(\G,(g_i)_{1 \leq i \leq \theta},(\chi_i)_{1 \leq i
\leq \theta})$ be a YD-datum.

\smallbreak We say that $\D$ is a {\em YD-datum of Cartan type}
$(a_{ij})$ if
\begin{align}
q_{ij}q_{ji} = q_{ii}^{a_{ij}},\; q_{ii} \neq 1,\; 0 \leq -a_{ij} < \ord(q_{ii}), {\rt \text{ for all } 1 \leq i,j \leq \theta,}
\end{align}
where the $q_{ij}$ are defined by \eqref{q}, and $1 \leq
\ord(q_{ii}) \leq \infty$.
\end{Def}

Note that the equivalence relation \eqref{sim} can be described as
usual in terms of the Cartan matrix.  For all $1 \leq i,j \leq
\theta$, $i \sim j$ if and only if there are vertices $i_1,\dots,
i_t \in \mathbb{I}, t \geq 2$ with $i_1 =i$, $i_t =j$,
$a_{i_l,i_{l+1}} \neq 0$ for all $1 \leq l < t$.

\begin{Def}\label{reducedCartanYD-datum}
A {\em reduced YD-datum of Cartan type}
$$\D(\G,(L_i)_{1\leq i \leq \theta},(K_i)_{1 \leq i \leq \theta},
(\chi_i)_{1 \leq i \leq \theta}, (a_{ij})_{1 \leq i,j \leq \theta})$$
is a reduced YD-datum $\D(\G,(L_i)_{1\leq i \leq \theta},(K_i)_{1
\leq i \leq \theta}, (\chi_i)_{1 \leq i \leq \theta})$ such that
for all $ 1 \leq i,j \leq \theta$
\begin{align}
q_{ij}q_{ji} &= q_{ii}^{a_{ij}},\quad q_{ii} \neq 1,\quad 0 \leq
-a_{ij} < \ord(q_{ii}),\label{Cartanreduced1}\end{align} where
$q_{ij} = \chi_j(K_i)$, as in page \pageref{qij}.
\end{Def}

We introduce an important condition which generalizes the notion of $X$-regular root data in \cite[Chapter 2]{L}.
\begin{Def}\label{Defregular}
A reduced YD-datum $$\D_{red}=\D(\G,(L_i)_{1\leq i \leq
\theta},(K_i)_{1 \leq i \leq \theta}, (\chi_i)_{1 \leq i \leq
\theta})$$ is called {\em regular}
 if the characters $\chi_1,\dots,\chi_{\theta}$ are $\mathbb{Z}$-linearly independent in $\widehat{\G}$.
\end{Def}

We fix a generic, see Definition \ref{YD-datum}, reduced YD-datum of
Cartan type
$$\D_{red}=\D_{red}(\G, (L_i)_{1 \leq i \leq \theta}, (K_i)_{1 \leq i \leq \theta}, (\chi_i)_{1 \leq i \leq \theta} ,(a_{ij})_{1 \leq i,j \leq \theta}).$$
 By \cite[Lemma
2.4]{AS-crelle} we can choose  $d_1,
\dots, d_{\theta}\in \N - 0$ such that
\begin{equation}\label{DJ1}
d_ia_{ij}=d_ja_{ji} \text{ for all } i,j \in \I.
\end{equation}
Let $\mathcal{X}$ be the set of connected components of $\I=\{1,\dots,\theta\}$ with respect to the Cartan matrix $(a_{ij})_{1 \leq i,j \leq \theta}$.

\smallbreak It is useful to single out the following subgroup of $\G$.

\begin{Def}\label{regular}
Let $\G^2$ be the subgroup of $\G$ generated by the products $K_1L_1,\dots,K_{\theta}L_{\theta}$.
\end{Def}

\begin{Lem}\label{independent}
\emph{(1)} Let $J \subset \I$ be a connected component. Then there
are $q_J \in \ku^{\times}$ which is not a root of unity, and roots
of unity $\zeta_j \in \ku$, $j \in J$, such that
$q_{jj}=q_J^{2d_j} \zeta_j$ for all $j \in J$. In particular, the
elements $(q_{jj})_{j \in J}$ are $\mathbb{N}$-linearly
independent, that is, if $(n_j)_{j \in J}$ is a family of natural
numbers, then $\prod_{j \in J} q_{jj}^{n_j} =1$ implies that $n_j
=0$ for all $j \in J$.

\medbreak \emph{(2)} If $(a_{ij})$ is invertible, \emph{e.~g.} if it is a Cartan matrix of finite
type, then $\D_{red}$ is regular.\label{finiteregular}

\medbreak \emph{(3)} If $\D_{red}$ is regular and the index of
$\G^2$ in $\G$ is finite, then the Cartan matrix $(a_{ij})$ is
invertible.

\medbreak \emph{(4)} {\rt If $(a_{ij})$ is a Cartan matrix of finite type, then for all connected components $J \subset \I$ there is an element $q_J \in \ku^{\times}$ such that}
{\rt \begin{align}
\rt q_{ii} &= q_J^{d_ia_{ij}} \text {for all }i \in J, \; J \in
\mathcal{X}.\label{finiteDJ}
\end{align}}
\end{Lem}
\pf (1) We choose an element $i \in J$, and $q_J \in \ku$ with
$q_{ii} = q_J^{2d_i}$. Then for all $j \in J$ there are $i_1,
\dots,i_t \in J$, $t \geq 2$, with $i_1=i$, $i_t=j$, and
$a_{i_li_{l+1}} \neq 0$ for all $1\leq l <t$. By applying
\eqref{Cartanreduced1} several times we obtain
$$q_{ii}^{a_{i_1i_2}a_{i_2i_3} \cdots a_{i_{t-1}i_t}} = q_{jj}^{a_{i_2i_1} a_{i_3i_2} \cdots a_{i_ti_{t-1}}}.$$
On the other hand
$$q_{ii}^{a_{i_1i_2}a_{i_2i_3} \cdots a_{i_{t-1}i_t}} = q_J^{2d_ia_{i_1i_2}a_{i_2i_3} \cdots a_{i_{t-1}i_t}}=q_J^{2d_ja_{i_2i_1}a_{i_3i_2} \cdots a_{i_ti_{t-1}}}$$
by applying \eqref{DJ1} several times. Hence for all $j \in J$ there is a root of unity $\zeta_j \in \ku$ such that $q_{jj}=q_J^{2d_j} \zeta_j$. In particular, the elements $(q_{jj})_{j \in J}$ are $\mathbb{N}$-linearly independent since $\D_{red}$ is generic, hence $q_J$ is not a root of unity.

(2) Suppose $n_1, \dots, n_{\theta}$ are integers with
$\prod_{i=1}^{\theta} \chi_i^{n_i} = 1.$ Let $J \subset \I$ be a
connected component. Since $\chi_i(K_jL_j) = {\rt q_{ii}}^{a_{ij}}$ for
all $i,j$ by \eqref{Cartanreduced1}, \eqref{reduced1} we obtain
for all $j \in J$
$$1=\prod_{i=1}^{\theta} \chi_i^{n_i}(K_jL_j) = \prod_{i \in J} q_{ii}^{n_ia_{ij}},$$
where in the last product we can assume that $i \in J$ since
$a_{ij}=0$ for all $i \notin J$. By the proof of (1) we may write
$q_{ii} = q_J^{2d_i} \zeta_i$ for all $ i \in J$, where the
$\zeta_i$ are roots of unity. Thus $\prod_{i \in J}
q_J^{2d_in_ia_{ij}} =1$ for all $j \in J$ since $q_J$ is not a
root of unity. {\rt Since $(a_{ij})_{i,j \in J}$ is invertible}, it follows that $n_i=0$ for all $i \in J$.

(3) Since $\chi_1,\dots,\chi_{\theta}$ are $\mathbb{Z}$-linearly
independent characters and $\G /\G^2$ is a finite group, the
restrictions of $\chi_1,\dots,\chi_{\theta}$ to the subgroup
$\G^2$ are $\mathbb{Z}$-linearly independent. Let $ J$ be any
connected component of $I$ with respect to $(a_{ij})$.  We use the
notation of the first part of the proof. Then $\chi_j(K_iL_i) =
q_{ii}^{a_{ij}} = q_J^{2d_i a_{ij}} \zeta_i$ for all $i,j \in J$.
Assume $n_j,j \in J$ are integers with $\sum_{j \in J} a_{ij}n_j
=0$ for all $i \in J$. Let $n\in \N$ with $\zeta_i^n=1$ for all $i
\in J$. Then
$$\prod_{j \in J}\chi_j^{nn_j}(K_iL_i) =
\prod_{j \in J} (q_J^{2d_ia_{ij}nn_j} \zeta_i^{nn_j}) =
q_J^{\sum_{j \in J} 2d_i a_{ij}nn_j} =1$$ for all $i \in J$. Since
the restrictions of the characters $\chi_j,j \in J$ to $\G^2$ are
$\mathbb{Z}$-linearly independent, and since $\chi_j(K_iL_i) =1$
for all $j \in J, i \in I \setminus J$, it follows that $n_j=0$
for all $j \in J$. Hence the matrix $(a_{ij})_{i,j \in J}$ is
invertible. Since $J$ was an arbitrary connected component, the
claim is proved.

{\rt (4) Finally, it is not difficult to see, by inspection,
that \eqref{finiteDJ} holds for data of finite Cartan type.} \epf

\begin{Rem}\label{DJ}
{\rt (1) The following relations hold in $\U = \Uc(\D_{red}, \ell)$ for all $1 \leq i,j \leq \theta, i \neq j$:
\begin{align}
\ad_l(E_i)^{1-a_{ij}}(E_j) &=\sum_{s=0}^{1 - a_{ij}} c_{ijs} E_i^s E_j E_i^{1 - a_{ij} - s}=0, \label{SerreE}\\
\ad_r(F_i)^{1-a_{ij}}(F_j) &=\sum_{s=0}^{1 - a_{ij}} d_{ijs} F_i^s F_j F_i^{1 - a_{ij} - s}=0, \label{SerreF}
\end{align}
where for all $1 \leq i,j \leq \theta, i \neq j, 0 \le s \le 1 - a_{ij}$, $c_{ijs},d_{ijs}$ are non-zero elements in $\ku$.
\pf The first equality in \eqref{SerreE} follows from the quantum binomial formula in $\End(\U^+)$, since for all $1 \le i \le \theta$, $\ad_l(E_i) = L_{E_i} - R_{E_i} \ad_lK_i$, and $(R_{E_i} \ad_lK_i) L_{E_i} = q_{ii} L_{E_i} (R_{E_i} \ad_lK_i)$, where $L_{E_i}(x) = L_i x$ and $R_{E_i}(x) = xE_i$ for all $x \in \U^+$. In same way the first equality in \eqref{SerreF} is shown. The second equality in \eqref{SerreE}  follows from \eqref{eqn:serre} since by definition the elements $v_i = E_i, 1 \leq i \leq \theta,$ satisfy the relations of the Nichols algebra $\NA(V)$.
By the same reason $\ad_l(w_i)^{1-a_{ij}}(w_j) =0$ for all $1 \leq i,j \leq \theta, i \neq j$. Hence $\Ss(\ad_l(w_i)^{1-a_{ij}}(w_j)) = \ad_r(\Ss(w_i))^{1-a_{ij}}(\Ss(w_j)) =0$ for all $1 \leq i,j \leq \theta, i \neq j$, where $\Ss$ is the antipode of the Hopf algebra $\NA(W) \# \ku [\Lambda]$. This proves the second equality in \eqref{SerreF} since by Corollary \ref{iota} (2) $\varphi^-(\Ss(w_i)) = - q_{ii} F_i$ for all $1 \le i \le \theta$.
\epf}
(2) Assume that the braiding matrix $(q_{ij})$ satisfies
\begin{align}
q_{ii} &= q_{J}^{2d_i} \quad \text{ for all }i \in J, \; J \in
\mathcal{X}\label{DJ2}.
\end{align}
Then $(q_{ij})$ is twist-equivalent to a braiding of
Drinfeld-Jimbo type \cite{AS-crelle}. Indeed, let
$\hat{q}_{ij}=q_J^{d_ia_{ij}}$, for all $J \in \mathcal{X}$ and
$i,j \in J$; set $\hat{q}_{ij}=1$, for $1 \leq i,j \leq \theta$
such that $i \nsim j$. Then $(\hat{q}_{ij})$ is of Drinfeld-Jimbo
type, and the braidings $(q_{ij})$ and $(\hat{q}_{ij})$ are
twist-equivalent since $q_{ij}q_{ji} = \hat{q}_{ij} \hat{q}_{ji},
q_{ii}=\hat{q}_{ii}$ for all $i,j$.

{\rt In this case, the braided Serre
relations \eqref{eqn:serre}
are defining relations of the Nichols algebras $\mathfrak{B}(V)$ and
$\mathfrak{B}(W)$. This follows by twisting from \cite[33.1.5]{L} when the elements $q_J \in \ku$ are transcendental, and from  \cite[Theorem 15]{Ro3} (see also \cite[Subsection 3.4]{HK}) when they  are not roots of unity.
Thus in Definition \ref{Unondegenerate} the relations \eqref{U1}, \eqref{U2} of $\U$ can be replaced by \eqref{SerreE} and \eqref{SerreF}.

To describe the relations  explicitly (cf. \cite[Lemma 1.6]{RS2}), let
\begin{equation}
p_{ij}= q_{ij}{\hat{q}_{ij}}^{\;-1}, \quad i \in J,\; J \in
\mathcal{X}, \;1 \leq j \leq \theta.
\end{equation}}
{\rt Then \eqref{SerreE}, \eqref{SerreF} are equivalent to
\begin{align} \sum_{s=0}^{1-a_{ij}}(-p_{ij})^s  \begin{bmatrix}
1-a_{ij}\\
s
\end{bmatrix}_{q_J^{d_i}} E_i^{1-a_{ij}-s} E_j E_i^s &=0, \label{explicitSerreE}\\
\sum_{s=0}^{1-a_{ij}}(-p_{ij})^s  \begin{bmatrix}
1-a_{ij}\\
s
\end{bmatrix}_{q_J^{d_i}} F_i^{s} F_j F_i^{1-a_{ij}-s} &=0,\label{explicitSerreF}
\end{align}
for all $i \in J$, $J \in \mathcal{X}$ and $1
\leq j \leq \theta$, $i \neq j$.}
\end{Rem}

\section{Representation theory of $\U$}\label{sect:$U$-modules}

In this section we assume that $\D_{red}$ is generic, regular, and
of Cartan type; we denote $\U = \Uc(\D_{red},\ell)$. {\rt We extend \cite[Sections  3.4 and 3.5]{L}.}

{\rt Let $Q$ be the subgroup of $\widehat{\G}$ generated by
$\chi_1,\dots,\chi_{\theta}$. Thus  by regularity $$\Z^\I
\xrightarrow{\cong} Q,\quad \alpha \mapsto \chi_{\alpha},$$ is bijective.}

\subsection{The category $\mathcal{C}^{hi}$}\label{subsection:C}

\

Let $\mathcal{C}$ be the full subcategory of $_{\U}\mathcal{M}$
consisting of all left $\U$-modules $M$ which are direct sums of
1-dimensional $\G$-modules, that is, which have a weight space
decomposition $M = \oplus_{\chi \in \widehat{\G}} M^{\chi}$, where
$$M^{\chi} = \{m \in M \mid gm=\chi(g)m \text{ for all } g \in \G \}$$
for all $\chi \in \widehat{\G}$. A character $\chi \in
\widehat{\G}$ is called a \emph{weight} for $M$ if $M^{\chi} \neq
0$.


\medbreak Let $\mathcal{C}^{hi}$ be the full subcategory of
$\mathcal{C}$ defined as follows. A module $M \in \mathcal{C}$ is
in $\mathcal{C}^{hi}$ if for any $m \in M$ there is an integer $N
\geq 0$ such that $\U^+_{\alpha}m =0$ for all $\alpha \in \N^{\I}$
with $|\alpha| \geq N$.

\medbreak Note that both categories $\mathcal{C}$ and
$\mathcal{C}^{hi}$ are closed under sub-objects and quotient
objects in $_{\U}\mathcal{M}$.

\medbreak We begin {\rt with } a technical result to be used later.

\begin{Prop}\label{corquasi}
Let $M \in \mathcal{C}^{hi}$. Then multiplication with $\Omega$ on
$M$  is a well-defined operator mapping each weight space of $M$
into itself. For all $\chi \in \widehat{\G},m \in M^{\chi}$, and
$1 \leq i \leq \theta$,
\begin{align*}
\Omega E_im&=(\chi\chi_i)(K_iL_i)^{-1} E_i\Omega m,\tag{1}\\
\Omega F_i m &= \chi(K_iL_i)F_i\Omega m.\tag{2}
\end{align*}
\end{Prop}
\pf For all $m \in M$, $\Omega m =\sum_{\alpha \in \N^{\I}}
\sum_{k=1}^{d_{\alpha}} S(x_{\alpha}^k) y_{\alpha}^k m$ is a
finite sum since $M \in \mathcal{C}^{hi}$. Hence multiplication
with $\Omega$ is a well-defined operator on $M$. For all $\alpha
\in \N^{\I}$ and $x \in \U^-_{-\alpha}, y \in \U^+_{\alpha}$ the
element $\Ss(x)y$ commutes with all $g \in \G$. Hence $\Omega : M
\to M$ is $\G$-linear and maps each weight space of $M$ into
itself.

To prove (1) let $\chi \in \widehat{\Gamma}$ and $m \in M^{\chi}$.
We apply $\Ss \o \id$ to Theorem \ref{quasi} (1), multiply and
obtain for all $\alpha \in \N^{\I}$
\begin{align*}
&\sum_{k=1}^{d_{\alpha}}\Ss(x_{\alpha}^k) S(E_i) y_{\alpha}^k +
\sum_{l=1}^{d_{\alpha - \alpha_i}} \Ss(x_{\alpha - \alpha_i}^l) K_i^{-1} E_i y_{\alpha - \alpha_i}^l \\
&= \sum_{k=1}^{d_{\alpha}} \Ss(E_i) \Ss(x_{\alpha}^k) y_{\alpha}^k
+ \sum_{l=1}^{d_{\alpha - \alpha_i}} L_i \Ss(x_{\alpha -
\alpha_i}^l)y_{\alpha - \alpha_i}^lE_i.
\end{align*}
Here both sums over $l$ are zero if $\alpha - \alpha_i \notin
\N^{\I}$. Since $\Ss(E_i) = -K_i^{-1}E_i$ it follows that

\begin{align*}
&-\sum_{\alpha \in \N^{\I}}\sum_{k=1}^{d_{\alpha}}
\Ss(x_{\alpha}^k)K_i^{-1} E_iy_{\alpha}^k m + \sum_{\alpha \in
\N^{\I}}
\sum_{l=1}^{d_{\alpha - \alpha_i}} \Ss(x_{\alpha - \alpha_i}^l) K_i^{-1} E_i y_{\alpha - \alpha_i}^l m\\
&= -\sum_{\alpha \in \N^{\I}}\sum_{k=1}^{d_{\alpha}} K_i^{-1}E_i
\Ss(x_{\alpha}^k) y_{\alpha}^k m + \sum_{\alpha \in
\N^{\I}}\sum_{l=1}^{d_{\alpha - \alpha_i}} L_i \Ss(x_{\alpha -
\alpha_i}^l)y_{\alpha - \alpha_i}^lE_im.
\end{align*}

Since the left hand side of the last equation is zero we obtain
\begin{align*}
0 &= -K_i^{-1} E_i \Omega m + L_i \Omega E_i m,\\
\intertext{hence}
\Omega E_i m &= (K_iL_i)^{-1} E_i\Omega m = (\chi\chi_i)(K_iL_i)^{-1} E_i \Omega m.
\end{align*}
In the same way (2) follows from Theorem \ref{quasi} (2).
\epf

Let $\chi \in \widehat{\G}$. We define the {\em Verma module}
$$M(\chi) = {\U}/ \big(\sum_{i=1}^{\theta} {\U}E_i + \sum_{g \in \G} {\U}(g - \chi(g)\big).$$
The inclusion $\U^- \subset \U$ defines a $\U^-$-module isomorphism
\begin{equation}\label{isoU^-}
\U^- \xrightarrow{\cong} M(\chi) = \U/\big(\sum_{i=1}^{\theta} \U
E_i + \sum_{g \in \G} \U (g-\chi(g))\big).
\end{equation}
This follows from the triangular decomposition in Corollary \ref{decomposition} (2).

\medbreak Let $m_{\chi}\in M(\chi)$ be the residue class of 1 in
$M(\chi)$. Then $M(\chi) \in \mathcal{C}$, $m_{\chi} \in
M(\chi)^{\chi}$, and $E_im_{\chi}=0$ for all $1 \leq i \leq
{\theta}$. The pair $(M(\chi),m_{\chi})$ has the following
universal property: For any $M \in \mathcal{C}$ with $m \in
M^{\chi}$ such that $E_im=0$ for all $1 \leq i \leq \theta$ there
exists a unique $\U$-linear map $t : M(\chi) \to M$ such that
$t(m_{\chi})=m$.

\medbreak The Verma module $M(\chi)$ and all its quotients belong
to the category $\mathcal{C}^{hi}$. We define a partial order $\leq$ on $\widehat{\G}$.

\begin{Def}
For all $\chi, \chi' \in \widehat{\G}$ we write $\chi' \leq \chi$ if there is an element $\alpha \in \N^{\I}$ such that $\chi = \chi' \chi_{\alpha}$.
\end{Def}
Note that $\leq $ is a partial order in $\widehat{\G}$ since $\D_{red}$ is regular.
\begin{Lem}\label{Verma}
Let $\chi \in \widehat{\G}$, and $M \in \mathcal{C}$.
Suppose $\chi$ is a maximal weight for $M$ and $m \in M^{\chi}$. Then $E_im =0$ for all $1 \leq i \leq \theta$, and $\U m$ is a quotient of $M(\chi)$.
\end{Lem}
\pf This follows from the universal property of the Verma module since $E_i m \in M^{\chi \chi_{\alpha_i}}$, and $M^{\chi \chi_{\alpha_i}} =0$ by maximality of $\chi$.
\epf

\medbreak Let $M \in \mathcal{C}$ and $C$ be a coset of $Q$ in $\widehat{\G}$. Then $M_{C} = \oplus_{\chi \in C} M^{\chi}$ is an object of $\mathcal{C}$. We note that $M=\oplus_{C} M_{C}$, where $C$ runs over the $Q$-cosets of $\widehat{\G}$.

By regularity, for all $\chi \in \widehat{\G}$
$$M(\chi) = \oplus_{\alpha \in \N^{\I}}M(\chi)^{\chi (\chi_{\alpha})^{-1}},\; M(\chi)^{\chi} = \ku m_{\chi},$$
since $M(\chi)$ is the $\ku$-span of the residue classes of $F_{i_1} \cdots F_{i_n}, 1 \leq i_1,\cdots, i_n \leq \theta, n \geq 0$. Thus $\chi$ is a weight of $M(\chi)$ with one-dimensional weight space, $\chi' \leq \chi$ for all weights $\chi'$ of $M(\chi)$, and $M(\chi) = (M(\chi))_C$, where $C=\chi Q$.
Because of these remarks, the proof of the following Lemma is standard.

\begin{Lem}\label{lem:simplehw}
If $\chi \in \widehat{\G}$, then $M(\chi)$ has a unique maximal submodule $M'(\chi)$; the quotient
$\Lambda(\chi) := M(\chi)/M'(\chi)$ is the unique (up to isomorphisms) simple module with highest weight $\chi$. \qed
\end{Lem}

\begin{Lem}\label{hi}
Suppose $M \in \mathcal{C}^{hi}$ is finitely generated as a $\U$-module.
\begin{enumerate}
\item The dimension of $M^{\chi}$ is finite for all $\chi \in \widehat{\G}$.
\item For all $\chi' \in \widehat{\G}$ there are only finitely many weights $\chi$ for $M$ which satisfy $\chi' \leq \chi$.
\item Every non-empty set of weights for $M$ has a maximal element.
\end{enumerate}
\end{Lem}
\pf We may assume that $M \neq 0$. In this case $M$ is generated
by weight vectors $v_1, \ldots , v_r$. Let $\chi_1, \ldots,\chi_r
\in \widehat{\Gamma}$ be the corresponding weights. Let
 $\chi$ be a weight for $M$. Observe that $M^{\chi}$ is
spanned
 by elements of the form
$$
0 \neq m = F_{i_1}\cdots F_{i_s}E_{j_1}\cdots E_{j_t}g\cdot v_i =
\chi_i(g) F_{i_1}\cdots F_{i_s}E_{j_1}\cdots E_{j_t}\cdot v_i,
$$
where $1 \leq i \leq r$, $g \in \Gamma$, $0 \leq s, t$,
$1 \leq i_p, j_q \leq \theta$ for all $1 \leq p \leq s$,
$1 \leq q \leq t$, and $\chi = \chi_{-\beta + \alpha} \chi_i$, where $\beta = \alpha_{i_1} + \cdots + \alpha_{i_s}$,
and  $\alpha = \alpha_{j_1} + \cdots + \alpha_{j_t}$.
Since $M \in {\mathcal
 C}^{hi}$
there are only finitely many $\alpha$'s for each $1 \leq i \leq r$, and for each pair $(\alpha,i)$ there is exactly one $\beta$ with  $\chi = \chi_{-\beta + \alpha} \chi_i$.
Here we use the fact that $\chi_1, \ldots,
\chi_\theta$ are $\mathbb{Z}$-linearly independent, that is for
 $\alpha,
 \beta
 \in \Z^\I$ the equations $\chi_\alpha = \chi_\beta$ implies
 $\alpha =
 \beta$. We have established (1).

 Let $\chi' \in \widehat{\G}$ and $\chi$ a weight for $M$ such that $\chi' \leq \chi$. Note that $ \chi \leq \chi_{\alpha} \chi_i$ for some $\alpha$ and $i$ as above. This proves (2) since there are only finitely many such pairs $(\alpha,i)$, and since for all $\chi_1,\chi_2 \in \widehat{\G}$ the segment $$[\chi_1,\chi_2]=\{\varphi \in \widehat{\G} \mid \chi_1 \leq \varphi \leq \chi_2 \}$$ is finite.
A consequence of (2) is that every chain of weights $\chi_1 \leq \chi_2
\leq \cdots$ is finite; hence (3) follows.
\epf

\subsection{Integrable modules}\label{subsection:integrable}

\

A left $\U$-module $M$ is called {\em integrable} if $M \in
\mathcal{C}$, and for any $m \in M$ and $1 \leq i \leq {\theta}$
there is a natural number $n \geq 1$ such that $E_i^nm = F_i^nm
=0$.


\medbreak The following notion from \cite{RS2} is an adaptation to
the present setting of the classical concept in Lie theory. A
character $\chi \in \widehat{\G}$ is called {\em dominant} if
there are natural numbers $m_i \geq 0$ such that $\chi(K_iL_i) =
q_{ii}^{m_i}$ for all $1 \leq i \leq {\theta}$.

\smallbreak We denote the set of all dominant characters in
$\widehat{\G}$ by $\widehat{\G}^+$.

\begin{Def}\label{def:dominant}
Let $\chi\in \widehat{\G}^+$ and $m_i \geq 0$ for all
$1 \leq i \leq \theta$ such that $\chi(K_iL_i) = q_{ii}^{m_i}$ for all $1 \leq i \leq \theta$.
Set
$$L_{\U}(\chi)= \U/\big(\sum_{i=1}^{\theta} \U E_i + \sum_{i=1}^{\theta} \U F_i^{m_i +1}
+ \sum_{g \in \G} \U (g - \chi(g))\big).$$
\end{Def}
We will write $L(\chi) = L_{\U}(\chi)$ when the Hopf algebra $\U$ is fixed.

\begin{Lem}\label{EF}
Let $n \geq 1$ and $1 \leq i \leq \theta$. Then
\begin{enumerate}
\item $E_iF_i^n = F_i^nE_i + \ell_i
\displaystyle{\frac{q_{ii}^n -1}{q_{ii} -1}} (K_i - L_i^{-1}
q_{ii}^{-n +1}) F_i^{n-1}.$\label{EF1}

\medbreak\item $F_iE_i^n = E_i^nF_i + \ell_i
\displaystyle{\frac{q_{ii}^n -1}{q_{ii} -1}} (L_i^{-1} - K_i
q_{ii}^{-n +1})E_i^{n-1}.$\label{EF2}

\medbreak\item $F_i^nF_j \in \sum_{s=0}^{-a{ij}} \ku F_i^s F_j
F_i^{n -s}$, if $n \geq 1 - a_{ij}$.\label{EF3}

\medbreak\item $E_i^nE_j \in \sum_{s=0}^{-a{ij}} \ku E_i^s E_j
E_i^{n -s}$, if $n \geq 1 - a_{ij}$.\label{EF4}
\end{enumerate}
\end{Lem}

\pf (1) and (2) follow from Prop. \ref{rulesE} \eqref{rule4} and
Prop. \ref{rulesF} \eqref{rule1}, or can be shown directly by
induction on $n$. (3) and (4) follow from the Serre relations
{\rt \eqref{SerreF} and \eqref{SerreE}} and the
observation that for $a,b$ in an algebra and $r \geq 1$ the
relation $a^rb \in \sum_{s=0}^{r-1} \ku a^sba^{r-s}$ implies $a^nb
\in \sum_{s=0}^{r-1} \ku a^sba^{n-s}$ for all $n \geq r$. \epf

Let $\chi$ be dominant, and let $\ell_{\chi}$ be the residue class
of 1 in $L(\chi)$. By the next lemma the pair
$(L(\chi),\ell_{\chi})$ has the universal property of the Verma
module with respect to integrable modules in $\mathcal{C}$.

\begin{Prop}\label{dominant}
Let $M \in \mathcal{C}$ be integrable and $\chi \in \widehat{\G}$.
Assume that there exists an element $0 \neq m \in M^{\chi}$ such
that $E_im=0$ for all $1 \leq i \leq \theta$. Then $\chi$ is
dominant, and there is a unique $\U$-linear map $t : L(\chi) \to
M$ such that $t(\ell_{\chi}) = m$.
\end{Prop}
\pf Let $1 \leq i \leq \theta$. Since $M$ is integrable, there is
an integer $n \geq 1$ such that $F_i^{n}m=0$, $F_i^{n-1}m \neq 0$.
By Lemma \ref{EF} (1)
\begin{align*}
0&= E_iF_i^nm = \ell_i
\displaystyle{\frac{q_{ii}^n -1}{q_{ii} -1}} (K_i - L_i^{-1} q_{ii}^{-n +1}) F_i^{n-1}m\\
&=\ell_i \displaystyle{\frac{q_{ii}^n -1}{q_{ii} -1}}
\left(\chi(K_i) q_{ii}^{-n +1} - \chi(L_i^{-1})\right) F_i^{n-1}m,
\end{align*}
since $F_i^{n-1}m \in M^{\chi \chi_i^{1-n}}$. Since $q_{ii}$ is
not a root of unity, it follows that $\chi(K_iL_i) =
q_{ii}^{n-1}$. Hence $n=m_i + 1$. Thus $\chi$ is dominant, and the
universal $\U$-linear map $M(\chi) \to M,\;m_{\chi} \mapsto m,$
factorizes over $L(\chi)$. \epf

\begin{Cor}\label{L} Let $\chi, \chi' \in \widehat{\G}^+$.

\emph{(1)} The isomorphism \eqref{isoU^-} induces an isomorphism
\begin{equation}\label{eqn:presentationL}
\U^-/\big(\sum_{i=1}^{\theta} \U^- F_i^{m_i+1} \big)
\xrightarrow{\cong} L(\chi),
\end{equation}
$L(\chi)$ is integrable, and $\dim L(\chi)^{\chi} = 1$, with basis
$\ell_{\chi}$.

\medbreak \emph{(2)}  The modules $L(\chi)$ and $L(\chi')$ are
isomorphic if and only if $\chi = \chi'$.
\end{Cor}

\pf (1) By Lemma \ref{dominant} $E_iF_i^{m_i +1} \overline{1} = 0$
in $M(\chi)$, $1 \leq i \leq \theta$. Hence the image of $\U
F_i^{m_i+1}$ in $M(\chi)$ coincides with the image of $\U^-
F_i^{m_i +1}$, and the map in (1) is bijective. In particular,
$L(\chi) \neq 0$, and $L(\chi)^{\chi}$ is one-dimensional with
basis $\overline{1}=l_{\chi}$. By Lemma \ref{EF} (3) $L(\chi)$ is
integrable. (2) follows from (1) since $L(\chi)$ and $L(\chi')$
have unique highest weights $\chi$ and $\chi'$. \epf

{\rt We note that in the proof of the last corollary we used the following rule in $\U^-$ to show that $L(\chi)$ is integrable: For all $1 \leq i,j \leq \theta, i \ne j,$ there are integers $n_{ij} \ge r_{ij} \geq 0$  such that
\begin{align}\label{weaker}
F^n_i F_j \in \U^- F_i^{n - r_{ij}} \; \text{ for all }n \geq n_{ij}.
\end{align}
This rule follows from Lemma \ref{EF} \eqref{EF3} with $r_{ij} = - a_{ij}, n_{ij} = 1 - a_{ij}$, that is, from the Serre relations which hold because $\D_{red}$ is of Cartan type. The assumption of Cartan type is only used here. Thus in Section \ref{sect:$U$-modules} we could replace it by \eqref{weaker}.}

\subsection{The quantum Casimir operator}\label{sectionCasimir}

\

We assume in this subsection {\rt the following condition on the diagonal entries of the braiding matrix $(q_{ij})$:
\begin{align}\label{Nli}
\text{If }\prod_{i=1}^{\theta} q_{ii}^{n_i} =1 , 0 \leq n_i  \in \mathbb{Z}, 1 \leq i \leq \theta,\text{ then } n_i = 0 \text{ for all }1 \leq i \leq \theta.
\end{align}
By Lemma \ref{independent} \eqref{Nli} holds} if $\I$ is
connected, that is if the Cartan matrix of $\D_{red}$ is
indecomposable. As in \cite[Chapter 6]{L} the next lemma is crucial for the
semisimplicity results.

\begin{Lem}\label{G}
Let $C$ be a coset of $Q$ in $\widehat{\G}$.
\begin{enumerate}
\item There is a function $G : C \to \ku^{\times}$ such that $G(\chi) = G(\chi \chi_i^{-1}) \chi(K_iL_i)$,
for all $ \chi \in C$ and $1 \leq i \leq \theta$.
$G$ is uniquely determined up to multiplication by a constant in
$\ku^{\times}$. \label{G1}

\item Let $G$ be as in \emph{(1)}. If $\chi,\chi' \in \widehat{\G}^+$ are
dominant characters with $\chi \geq \chi'$ and $G(\chi) = G(\chi')$,
then $\chi = \chi'$.\label{G2}
\end{enumerate}
\end{Lem}

\pf \eqref{G1} Let $C=\bar{\chi}Q$ where $\bar{\chi}$ is a fixed
element in the coset $C$, and pick $G(\bar{\chi})\in
\ku^{\times}$. For all $\alpha=\sum_{i=1}^{\theta} n_i \alpha_i
\in \Z^\I$ we define
\begin{align}
q_{\alpha}&= \prod_{i=1}^{\theta} q_{ii}^{n_i(n_i+1)} \prod_{1 \leq i < j \leq \theta} (q_{ij}q_{ji})^{n_in_j},\label{defq}\\
G(\bar{\chi}\chi_{\alpha}) &= G(\bar{\chi})\bar{\chi}(K_{\alpha}L_{\alpha})q_{\alpha}.\label{defG}
\end{align}
We first show that for all $\alpha \in \Z^\I, 1 \leq p \leq \theta$,
\begin{align}
q_{\alpha} = q_{\alpha - \alpha_p} \chi_{\alpha}(K_pL_p).\label{formulaq}
\end{align}
By definition
\begin{align*}
q_{\alpha - \alpha_p}&= \prod_{\substack{1 \leq i \leq \theta,\\i \neq p}}q_{ii}^{n_i(n_i+1)} q_{pp}^{(n_p-1)n_p} \prod_{\substack{1 \leq i < j \leq \theta,\\ i \neq p,j \neq p}} (q_{ij}q_{ji})^{n_in_j}
\prod_{\substack{1 \leq i \leq \theta,\\ i \neq p}} (q_{ip}q_{pi})^{n_pn_i - n_i}\\
&=\prod_{1 \leq i \leq \theta}q_{ii}^{n_i(n_i+1)} q_{pp}^{-2n_p}\prod_{1 \leq i < j \leq \theta} (q_{ij}q_{ji})^{n_in_j} \prod_{\substack{1 \leq i \leq \theta,\\ i \neq p}}(q_{ip}q_{pi})^{-n_i}\\
&=q_{\alpha} \chi_{\alpha}(K_pL_p)^{-1},
\end{align*}

\noindent where the last equality follows from $\chi_i(K_pL_p) =
q_{ip}q_{pi}$ for all $1 \leq i \leq \theta$.

\smallbreak It follows from \eqref{formulaq} that the function $G$
defined by \eqref{defG} has the desired property since for all
$\alpha \in \Z^\I, 1 \leq p \leq \theta$,
\begin{align*}
G(\bar{\chi}\chi_{\alpha}\chi_p^{-1})
(\bar{\chi}\chi_{\alpha})&(K_pL_p) =G(\bar{\chi})\bar{\chi}
(K_{\alpha-\alpha_p}L_{\alpha-\alpha_p}) q_{\alpha-\alpha_p}\bar{\chi}(K_pL_p) \chi_{\alpha}(K_pL_p)\\
&=G(\bar{\chi})\bar{\chi}(K_{\alpha}L_{\alpha})\bar{\chi}(K_p^{-1}L_p^{-1})
q_{\alpha - \alpha_p}
\bar{\chi}(K_pL_p) \chi_{\alpha}(K_pL_p)\\
&= G(\bar{\chi})\bar{\chi}(K_{\alpha}L_{\alpha})q_{\alpha-\alpha_p}\chi_{\alpha}(K_pL_p)\\
&=G(\bar{\chi})\bar{\chi}(K_{\alpha}L_{\alpha})q_{\alpha}\\
&= G(\bar{\chi}\chi_{\alpha}).
\end{align*}
The functions $G$ in \eqref{G1} are clearly unique up to a
non-zero scalar.

\smallbreak \eqref{G2} (a) We show by induction on $n \geq 0$ that
for all $1 \leq i_1,i_2,\dots,i_n \leq \theta$,

\begin{align}
\frac{G(\chi)}{G(\chi \chi_{i_1}^{-1} \chi_{i_2}^{-1} \cdots
\chi_{i_n}^{-1})} = \prod_{1 \leq p \leq n} \chi(K_{i_p}L_{i_p})
\prod_{1 \leq p<q \leq n}
\chi_{i_p}(K_{i_q}L_{i_q})^{-1}.\label{Ga}
\end{align}

This is clear for $n=0$ and follows by induction and \eqref{G1}
from
\begin{align*}
&G(\chi)G(\chi \chi_{i_1}^{-1} \cdots \chi_{i_{n+1}}^{-1})^{-1}\\
&=G(\chi)G(\chi \chi_{i_1}^{-1} \cdots \chi_{i_{n}}^{-1})^{-1}G(\chi \chi_{i_1}^{-1} \cdots \chi_{i_{n}}^{-1})G(\chi \chi_{i_1}^{-1} \cdots \chi_{i_{n+1}}^{-1})^{-1}\\
&=\prod_{1 \leq p \leq n} \chi(K_{i_p}L_{i_p}) \prod_{1 \leq p<q \leq n} \chi_{i_p}(K_{i_q}L_{i_q})^{-1}(\chi \chi_{i_1}^{-1} \cdots \chi_{i_n}^{-1})(K_{i_{n+1}}L_{i_{n+1}})\\
&=\prod_{1 \leq p \leq n+1} \chi(K_{i_p}L_{i_p})\prod_{1 \leq p<q \leq n+1} \chi_{i_p}(K_{i_q}L_{i_q})^{-1}.
\end{align*}

\smallbreak (b) Now we prove (2). By assumption there are indices
$1 \leq i_1,i_2,\dots,i_n \leq \theta$, $n \geq 0$, with
$\chi'=\chi\chi_{i_1}^{-1} \cdots \chi_{i_n}^{-1}$. Since $\chi$
and $\chi'$ are both dominant there are natural numbers $m_i,m'_i
\geq 0$, $1 \leq i \leq \theta$, such that
\begin{align}
&\chi(K_iL_i) = q_{ii}^{m_i},\;\chi'(K_iL_i) = q_{ii}^{m'_i} \text{ for all } 1 \leq i \leq \theta.\label{G3}
\end{align}
By assumption $G(\chi) = G(\chi')=G(\chi\chi_{i_1}^{-1} \cdots \chi_{i_n}^{-1})$. Hence by (a),
\begin{align}
&\prod_{1 \leq p \leq n} \chi(K_{i_p}L_{i_p}) = \prod_{1 \leq p<q \leq \theta} \chi_{i_p}(K_{i_q}L_{i_q}).\label{G4}
\end{align}
Then we obtain
\begin{align*}
\prod_{1 \leq p \leq n}q_{i_pi_p}^{m_{i_p} + m'_{i_p}}=&
\prod_{1 \leq p \leq n} \chi(K_{i_p}L_{i_p})\prod_{1 \leq p \leq n} \chi'(K_{i_p}L_{i_p})\\
=&\prod_{1 \leq p \leq n} \chi(K_{i_p}L_{i_p})^2 \prod_{1 \leq p,q \leq n} \chi_{i_q}^{-1}(K_{i_p}L_{i_p})\\
=\prod_{1 \leq p<q \leq n}\chi_{i_p}(K_{i_q}L_{i_q})^2 &\prod_{1
\leq p<q \leq n} \chi_{i_p}(K_{i_q}L_{i_q})^{-2}
\times \prod_{1 \leq p \leq n} \chi_{i_p}^{-1}(K_{i_p}L_{i_p})\\
=& \prod_{1 \leq p \leq n} q_{i_pi_p}^{-2}.
\end{align*}

\noindent For the third equality we used \eqref{G4} and that
$\chi_i(K_jL_j)= q_{ij}q_{ji} = \chi_j(K_iL_i)$ for all $1\leq i,j
\leq \theta$. By {\rt \eqref{Nli}} the family $(q_{ii})_{1 \leq i \leq
\theta}$ is $\mathbb{N}$-linearly independent and we get a
contradiction except $n=0$, that is $\chi = \chi'$. \epf

\begin{Ex}
Let $\G=\langle K_1,K_2\rangle$ be a free abelian group with basis $K_1$,
$K_2$, and $ 0 \neq q \in \ku$ not a root of unity. Let $L_1=K_1$,
$L_2=K_2$, and define characters $\chi_1,\chi_2 \in \widehat{\G}$
by
$$\chi_1(K_1)=q,\quad \chi_1(K_2)=1,\quad \chi_2(K_1)=1,\quad \chi_2(K_2)=q^{-1}.$$
Thus $\D_{red} =\D_{red}(\G, (L_i),(K_i), (\chi_i),(a_{ij}))$ is a
generic reduced YD-datum of Cartan type where
$a_{11}=a_{22}=2,a_{12}=a_{21}=0$, and $q_{11}q_{22}=1$. Define
$\chi, \chi' \in \widehat{{\rt \G}}$ by $\chi'(K_1)=q$,
$\chi'(K_2)=q^{-1}$, and $\chi = \chi' \chi_1 \chi_2$. Then $\chi'
\leq \chi$, and both are dominant.
{\rt Let $G$ be a function satisfying Lemma \ref{G} (1) for the coset $C=\chi'Q$.
Then
$$G(\chi)=G(\chi'\chi_1\chi_2) = G(\chi') (\chi'\chi_1)(K_1L_1) (\chi' \chi_1\chi_2)(K_2L_2) = G(\chi').$$} Thus Lemma \ref{G} (2) does not hold without
the assumption that the $q_{ii}$'s are $\mathbb{N}$-linearly
independent.
\end{Ex}

\begin{Prop}\label{OG}
Let $C$ be a coset of $Q$ in $\widehat{\G}$, and $M \in
\mathcal{C}^{hi}$ such that $M=M_C$. Choose a function $G$ as in
Lemma \ref{G} and define a $\ku$-linear map $\Omega_G : M \to M$
by $\Omega_G(m) = G(\chi)\Omega(m)$ for all $m \in M^{\chi},\chi
\in C$.
\begin{enumerate}
\medbreak\item The map $\Omega_G$ is $\U$-linear and locally
finite.

\medbreak\item If $0 \neq m \in M$ generates a quotient of a Verma
module $M(\chi)$ for some $\chi \in \widehat{\G}$, then $\chi \in
C$, and $\Omega_G(m) = G(\chi)m$.

\medbreak\item The eigenvalues of $\Omega_G$ are the $G(\chi)'s$,
where $\chi$ runs over the maximal weights of the submodules $N$
of $M$, in which case $\Omega_G (n) =G(\chi)n$ for all $n \in
N^\chi$.
\end{enumerate}
\end{Prop}
\pf (1) By Proposition \ref{corquasi} $\Omega_G$ is well-defined and
maps each weight space of $M$ to itself. Hence $\Omega_G : M \to
M$ is $\G$-linear. Let $1 \leq i \leq \theta, \chi \in
\widehat{\G},$ and $m \in M^{\chi}$. By Propositon \ref{corquasi}
\begin{align*}
\Omega_G(E_im)&=G({\chi_i}\chi)\Omega(E_im)= G(\chi_i\chi)(\chi\chi_i)^{-1}(K_iL_i) E_i\Omega m,\\
\Omega_G(F_im) &=G(\chi_i^{-1}\chi) \Omega(F_im) = G(\chi_i^{-1}\chi) \chi(K_iL_i) F_i\Omega m.
\end{align*}
On the other hand, $E_i\Omega_G(m) = G(\chi)E_i\Omega m,
F_i\Omega_G(m)= G(\chi)F_i\Omega m$. By Lemma \ref{G} (1),
$$G(\chi) =G(\chi\chi_i^{-1}) \chi(K_iL_i), G(\chi\chi_i) = G(\chi)
(\chi\chi_i)(K_iL_i).$$ Hence it follows that
$$\Omega_G(E_im)= E_i \Omega_G(m),\; \Omega_G(F_im)= F_i \Omega_G(m),$$
and we have shown that $\Omega_G$ is $\U$-linear. We show that $M$
is the sum of finite-dimensional $\Omega_G$-invariant subspaces.
Since any $\U$-submodule of $M$ is $\Omega_G$-invariant we may
assume that $M$ is finite-dimensional. In this case $M$ is the sum
of finite-dimensional weight spaces by Lemma \ref{hi} (1), and
weight spaces are $\Omega_G$-invariant. Hence $\Omega_{G}$ is a
locally finite linear map.

\medbreak (2) Write $\U \cdot m = \U \cdot n$, where $n \in
M^\chi$ and $E_i \cdot n = 0$ for all $1 \leq i \leq \theta$. Then
$\Omega_G (n) = G(\chi)n$ by definition of $\Omega_G$ and
consequently $\Omega_G(n') = G(\chi)n'$ for all $n' \in \U\cdot n$
since eigenspaces of module endomorphisms are submodules.

\medbreak (3) First of all, $G(\chi)$ is an eigenvalue of
$\Omega_G$ when $\chi$ is a maximal weight of submodule of $M$ by
Lemma 5.3 and part (2).

Conversely, suppose that $\lambda$ is an eigenvalue of $\Omega_G$
and $0 \neq m \in M$ satisfies $\Omega_G(m)  = \lambda m$. Since
$N = U \cdot m \neq 0$ is finitely generated, and $N \in {\mathcal
C}^{hi}$, by Lemma 5.4 (3) there is a maximal weight $\chi$ for
$N$. By Lemma 5.3 and part (2) we conclude that $G(\chi)$ is an
eigenvalue for the restriction $\Omega_G | N$. Since the
eigenvectors for $\Omega_G$ belonging  to
 $\lambda$ form a submodule of $M$, $G(\chi) = \lambda$.
\epf

The function $\Omega_G : M \to M$ in Proposition \ref{OG} is called the
{\em quantum Casimir operator}.

\subsection{Irreducible highest weight modules}\label{subsect:irrep}

\

\begin{Lem}\label{components}
Let $\chi \in \widehat{\G}^+$. Let $J$ be a connected component of
$\I$, $\I'=\I \setminus J$, and let $\U_J$ be the subalgebra of
$\U$ generated by $\G$ and $E_j,F_j, j \in J$ and $\U'$  the
subalgebra of $\U$ generated by $\G$ and $E_i,F_i, i \in \I'$.
Then the map
\begin{equation}\label{eqn:components}
\Phi : L_{\U_J}(\chi) \otimes L_{\U'}(\chi) \to L(\chi),\quad
\overline{u} \otimes \overline{u'} \mapsto \overline{uu'},
\end{equation} for
all $u \in \U^-_J$, $u' \in \U'^-$, is a $\ku$-linear isomorphism;
and
\begin{enumerate}
\medbreak\item $\Phi(g m \otimes g m') = \chi(g)\, g \Phi(m
\otimes m')$,

\medbreak\item $\Phi(E_j m \otimes m') = (\chi \psi^{-1})(K_j)\,
E_j \Phi(m \otimes m')$, if $m' \in L_{\U'}(\chi)^{\psi}$, $\psi
\in \widehat{\G}$,

\medbreak\item $\Phi(m \otimes E_i m') = E_i \Phi(m \otimes m')$,
\end{enumerate}
for all $j \in J$, $i \in \I'$, $m \in L_{\U_J}(\chi)$, $m' \in
L_{\U'}(\chi)$.
\end{Lem}

\pf The multiplication map defines an isomorphism $\U_J^- \otimes
\U'^- \to \U^-$ since the generators $F_j, j \in J$, of $\U_J^-$
and $F_i, i \in \I'$, of $\U'^-$ skew-commute. The kernel of the
canonical map
$$\U_J^- \otimes \U' \to \U_J^-/ \big(\sum_{j \in J} \U_J^- F_j^{m_{j}+1}\big)
\otimes \U'^-/\big(\sum_{i \in \I'} \U'^- F_i^{m_{i}+1}\big)$$ has
image
$$\sum_{j \in J} \U_J^-F_j^{m_j +1} \U'
+ \U_J^- \sum_{i \in \I'} \U'^- F_i^{m_i +1} = \sum_{i \in \I} \U F_i^{m_i +1}$$
under the multiplication map. Hence the induced map
$$\U_J^-/\big(\sum_{j \in J} \U_J^- F_j^{m_{j}+1}\big)
\otimes \U'^-/\big(\sum_{i \in \I'} \U'^- F_i^{m_{i}+1}\big) \to
\U^-/(\sum_{i \in \I} \U^-F_i^{m_{i}+1})$$ is bijective. Then
\eqref{eqn:components} is an isomorphism of vector spaces, by Corollary
\ref{L}.

\medbreak To prove (1) -- (3), we may assume that $m= u \,
\overline{1}$ and $m' = u' \, \overline{1}$, where $u \in
(\U_J^-)_{-\alpha}$, $u' \in (\U'^-)_{- \beta}$ are homogeneous
 with $\alpha \in \N^{J}$, $\beta \in
\N^{\I'}$.

\medbreak (1) Let $g \in \G$. Then
\begin{align*}
\Phi(gu\overline{1} \otimes gu'\overline{1}))
 &= (\chi_{-\alpha}\chi)(g) (\chi_{-\beta}\chi)(g) \Phi(u\overline{1} \otimes  u'\overline{1}) \\
&=(\chi_{-\alpha}\chi)(g) (\chi_{-\beta}\chi)(g) u u'\overline{1}\\
&= \chi(g)g u u'\overline{1}\\
&= \chi(g)g\Phi(u\overline{1} \otimes  u'\overline{1}).
\end{align*}

(2) Let $j \in J$.
We first note that there is an element $\widetilde{u} \in \U_J$ which is a $\ku$-linear combination of monomials in $F_l, l \in J$ and $K_j - L_j^{-1}$ where  in each monomial the factor $K_j - L_j^{-1}$ occurs exactly one, and such that
$E_j u =u E_j + \widetilde{u}.$
This follows by induction on $|\alpha|$, since
$$E_jF_ku = F_kE_ju + \delta_{jk}l(K_j - L_j^{-1})u=F_kuE_j +F_k\widetilde{u} + \delta_{jk}l(K_j - L_j^{-1})u$$
for all $k \in J$ by induction and \eqref{reducedlinking}. Then
\begin{align*}
\Phi(E_j u \overline{1} \otimes u' \overline{1}) &= \Phi(\widetilde{u} \overline{1} \otimes  u' \overline{1}),\\
E_j \Phi(u \overline{1} \otimes u' \overline{1}) &= E_j u u' \overline{1} =(u E_j + \widetilde{u}) u' \overline{1} = \widetilde{u} u' \overline{1},
\end{align*}
since $E_j  \overline{1} =0$ in $L_{\U_J}(\chi)$, and $E_j u' \overline{1} = u' E_j \overline{1} =0$ in $L_{\U'}(\chi)$. Hence (2) is equivalent to
\begin{equation}\label{2}
\Phi(\widetilde{u} \overline{1} \otimes  u' \overline{1}) = \chi_{\beta}(K_j) \widetilde{u} u' \overline{1},
\end{equation}
since $u'\overline{1} \in L_{\U'}(\chi)^{\chi \chi_{-\beta}}$,
hence $\chi \psi^{-1} = \chi_{\beta}$. To prove \eqref{2} we may
assume that $\widetilde{u} =u_1(K_j - L_j^{-1})u_2$, where $u_1
\in \U_j^-$ and $u_2 \in (\U_J^-)_{-\gamma},\gamma \in
\N^{J}$. Then
$$\widetilde{u}\overline{1} =u_1(K_j - L_j^{-1})u_2 \overline{1}
= ((\chi \chi_{- \gamma})(K_j) -(\chi \chi_{- \gamma})(L_j)^{-1})
u_1u_2 \overline{1}$$ in $L_J(\chi)$, and $\Phi(\widetilde{u}
\overline{1} \otimes  u' \overline{1})=((\chi \chi_{-
\gamma})(K_j) -(\chi \chi_{- \gamma})(L_j)^{-1}) u_1u_2u'
\overline{1}$.

Since $\chi_i(K_j) = q_{ji} = q_{ij}^{-1}$, and $\chi_i(L_j)^{-1} = \chi_j(K_i)^{-1} = q_{ij}^{-1}$ for all $i \in \I'$ it follows that in $L(\chi)$
\begin{align*}
\widetilde{u} u' \overline{1}&= u_1(K_j - L_j^{-1})u_2 u'\overline{1}\\
&= ((\chi \chi_{-\gamma} \chi_{-\beta})(K_j) -(\chi \chi_{-\gamma} \chi_{-\beta})(L_j)^{-1}) u_1 u_2 u' \overline{1} \\
&=\chi_{-\beta}(K_j)((\chi \chi_{-\gamma})(K_j) -(\chi \chi_{-\gamma})(L_j)^{-1}) u_1 u_2 u' \overline{1}\\
&=\chi_{-\beta}(K_j) \Phi(\widetilde{u} \overline{1} \otimes  u' \overline{1}).
\end{align*}

(3) As in the proof of (2) let $\widetilde{u'} \in \U'$ with $E_iu = u E_i + \widetilde{u'}$.
Then
\begin{align*}
\Phi(u \overline{1} \otimes E_i u' \overline{1}) &= \Phi(u \overline{1} \otimes \widetilde{u'} \overline{1}),\\
E_i \Phi(u \overline{1} \otimes u' \overline{1}) &= E_i uu' \overline{1} = uE_iu' \overline{1} =u \widetilde{u'} \overline{1}.
\end{align*}
To prove that $\Phi(u \overline{1} \otimes \widetilde{u'})= u \widetilde{u'} \overline{1}$ we may  assume that
$$\widetilde{u'} = u_1' (K_i - L_i^{-1})u_2', u_2' \in (\U'^-)_{-\delta}, \delta \in \N^{\I'}.$$
Then $\Phi(u \overline{1} \otimes \widetilde{u'} \overline{1})= ((\chi \chi_{-\delta})(K_i) -(\chi \chi_{-\delta})(L_i)^{-1}) u u_1u_2 \overline{1} =u \widetilde{u'} \overline{1}$.
\epf

The following theorem shows that $L(\chi)$  for $\chi \in \widehat{\G}^+$ coincides with the
simple module $\Lambda(\chi)$ from Lemma \ref{lem:simplehw}, thus
providing the defining relations of $\Lambda(\chi)$; it also
implies that the weight multiplicities of $\Lambda(\chi)$ are as
in the classical case for data of finite Cartan type.

\begin{Rem}\label{Kac}
 {\rt Let $(a_{ij})_{1 \le i,j \le \theta}$ be a symmetrizable Cartan matrix, $(\mathfrak{h}, \Pi,\Pi^{\vee})$ a realization of  $(a_{ij})_{1 \le i,j \le \theta}$ with  $\Pi=\{\alpha_1,\dots,\alpha_{\theta}\},\Pi^{\vee}=\{\alpha_1^{\vee},\dots,\alpha_{\theta}^{\vee}\}$ and $\mathfrak{g}$  the corresponding Kac-Moody Lie algebra (see \cite{K}). let $0 \ne q \in \ku$ be not a root of unity and $(V,c)$  a braided vector space with basis $v_1,\dots,v_{\theta}$ and braiding $c(v_i \ot v_j) = q^{d_ia_{ij}} v_j \ot v_i$ for all $1 \le i,j \le \theta$, where $(d_ia_{ij})$ is the symmetrized Cartan matrix.

Let $0 \le m_1,\dots,m_{\theta} \in \mathbb{Z}$. Choose $\lambda \in \mathfrak{h}^*$ with $\lambda(\alpha_i^{\vee}) = m_i$ for all $1 \le i \leq \theta$. Thus $\lambda$ is an integral weight of $\mathfrak{g}$. Let  $L(\lambda)$ be the irreducible $\mathfrak{g}$-module with}  highest weight $\lambda$. Then the multiplicities of the weight spaces of $L(\lambda)$ are given by
$$L(\lambda)_{\lambda - \alpha} \cong (U(\mathfrak{n}_+)/(\sum_{i=1}^{\theta} U(\mathfrak{n}_+) e_i^{m_i+1})_{\alpha}
\cong (\NA(V) /(\sum_{i=1}^{\theta} \NA(V) v_i^{m_i +1}))_{\alpha},$$
{\rt where $\deg(e_i) = \deg(v_i)= \alpha_i$ for all $i$, and $\alpha= \sum_{i=1}^{\theta} n_i \alpha_i, 0 \le n_i \in \mathbb{Z}$ for all $i$.
The first isomorphism is \cite[10.4.6]{K}, and the second isomorphism follows from \cite[33.1.3]{L} if $q$ is transcendental, and can be derived from \cite[Section 3.4]{HK} if $q$ is not a root of unity.}
\end{Rem}

\begin{Theorem}\label{mainsimple}
Let $\chi \in \widehat{\G}^+$.
\begin{enumerate}
\item $L(\chi)$ is a simple $\U$-module.

\item Any weight vector of $L(\chi)$ which is annihilated by all
$E_i, 1\leq i \leq \theta,$ is a scalar multiple of $\ell_{\chi}$.

\item If $(q_{ij})$ satisfies \eqref{DJ2}, in particular if the
Cartan matrix is of finite type, then the weight multiplicities
are as in the classical case, that is, given by the Weyl-Kac
character formula.
\end{enumerate}
\end{Theorem}

\pf We proceed by induction on the number of connected components
of $\I$. We first assume that $\I$ is connected. Then the results
of Subsection \ref{sectionCasimir} apply. Recall that $L(\chi) =
L(\chi)_C$ for the coset $C=\chi Q$.

\medbreak (1) Let $M$ be a non-zero submodule of $L(\chi)$. By
Lemma \ref{hi} (3) there is a maximal weight for $M$ since
$L(\chi)$ is finitely generated. Let $\chi'$ be such a weight.
Then $G(\chi')$ is an eigenvalue for $\Omega_G$ by Prop. \ref{OG}
(3). By part (2) of the same $G(\chi) = G(\chi')$. Since $L(\chi)$
is integrable $M$ is also. By Lemma \ref{Verma} and Prop.
\ref{dominant} $\chi'$ is dominant. Thus $\chi = \chi'$ by Lemma
\ref{G} (2). Since $L(\chi)^\chi$ is one-dimensional $L(\chi)^\chi
= M^{\chi}$. Thus $M = L(\chi)$ since $L(\chi)^\chi$ generates
$L(\chi)$. Thus we have shown that $L(\chi)$ is simple.

\medbreak (2) Let $\chi' \in \widehat{\G}$ and $0 \neq m \in
L(\chi)^{\chi'}$ such that $E_im=0$ for all $1 \leq i \leq
\theta$. By Proposition \ref{dominant} $\chi'$ is dominant and
there is a  $\U$-linear map $L(\chi') \to L(\chi)$ mapping
$\ell_{\chi'}$ onto $m$. This map is an isomorphism since
$L(\chi')$ and $L(\chi)$ are simple by the first part of the
proof. Hence $\chi' = \chi$ by Corollary \ref{L} (2), and $m$ is a
scalar multiple of $\ell_{\chi}$ by Corollary \ref{L} (1).

\medbreak {\rt (3) Let $\chi \in \widehat{\G}^+$ and $\chi(K_iL_i) = q_{ii}^{m_i}, 0 \le m_i \in \mathbb{Z}$, for all $1 \le i \le \theta$. By Corollary \ref{iota} the weights of $L(\chi)$ have the form $\chi \chi_{-\alpha}, \alpha \in \mathbb{N}^\I$, where for all $\alpha \in \mathbb{N}^\I$,
$L(\chi)^{\chi \chi_{-\alpha}}\cong \big(\U^-/\big(\sum_{i=1}^{\theta} \U^- F_i^{m_i+1}\big) \big)_{-\alpha}.$}
{\rt The bijective map $\widetilde{\kappa} : \NA(W) \to \U^-$ in Corollary \ref{iota} (3) induces an isomorphism
$\big(\NA(W)/\big(\sum_{i=1}^{\theta} \NA(W)w_i^{m_i + 1}\big)\big)_{\alpha} \cong \big(\U^-/\big(\sum_{i=1}^{\theta} \U^- F_i^{m_i+1}\big) \big)_{-\alpha}$ for all $\alpha \in \mathbb{N}^\I$
if we define $\deg(w_i) = \alpha_i$ for all $i$.

By assumption on the braiding (and since $\I$ is connected), $q_{ii} = q^{2d_i}$ for all $1 \le i \le \theta$, where $ 0 \ne q \in \ku$ is not a root of unity. The braiding matrix
$(q_{ji}^{-1})_{1 \leq i,j \leq \theta}$ of $W$ with respect to the basis $w_1,\dots,w_{\theta}$ is  twist equivalent to $(q^{-d_ia_{ij}})$. Let $\widehat{W}$ be the braided vector space with braiding matrix $(q^{-d_ia_{ij}})$ with respect to a basis $\widehat{w_1},\dots,\widehat{w_{\theta}}$. By \cite[Proposition 3.9, Remarks 3.10]{AS-cambr} $\NA(W)/\big(\sum_{i=1}^{\theta}\NA(W)w_i^{m_i + 1}\big) \cong \NA(\widehat{W})/\big(\sum_{i=1}^{\theta}\NA(\widehat{W})\widehat{w_i}^{m_i + 1}\big)$ as $\mathbb{Z}^\I$-graded vector spaces, where $\deg(\widehat{w_i}) = \alpha_i$ for all $i$.
The claim now follows from Remark \ref{Kac}.}

\medbreak

Now let $J$ be a connected component of $\I$ and $\I'=\I
\setminus J$. Let $\U_J$ be the subalgebra of $\U$ generated by
$\G$ and $E_j,F_j, j \in J$. Let $\U'$ be the subalgebra of $\U$
generated by $\G$ and $E_i,F_i, i \in \I'$. We assume by induction
that $L_{\U'}(\chi)$ satisfies (1), (2) and (3).

\medbreak We first show that $L(\chi)$ satisfies (2). Let $m \in
L(\chi)$ be a weight vector of weight $\chi' \in \widehat{\G}$
such that $E_i m =0$ for all $i \in \I$. By Lemma \ref{components}
(1) $\Phi$ induces a linear isomorphism of weight spaces
$$\bigoplus_{\substack{\chi^{-1} \varphi \psi = \chi'\\
\varphi, \psi \in \widehat{\G}}} L_{\U_J}(\chi)^{\varphi} \otimes L_{\U'}(\chi)^{\psi} \to L(\chi)^{\chi'}.$$

\noindent Hence there are finitely many elements $m_l \in
L_{\U_J}(\chi)^{\varphi_l}$, and $m'_l \in
L_{\U'}(\chi)^{\psi_l}$, $1 \leq l \leq n$,  with
$\varphi_l,\psi_l \in \widehat{\G}, \varphi_l \psi_l = \chi \psi$
for all $1 \leq l \leq n$ such that $m'_1,\dots,m'_n$ are
$\ku$-linearly independent and $\Phi(\sum_{l=1}^n m_l \otimes
m'_l) = m$. By Lemma \ref{components} (2)
$$\Phi\big(\sum_{l=1}^n (\chi_{-1}\psi_l)(K_j)E_jm_l \otimes m'_l \big)
= \sum_{l=1}^n E_j \Phi(m_l \otimes m_l') = E_j m =0.$$ for all $j
\in J$. Since $\Phi$ is bijective, and the elements $m'_l$ are
linearly independent it follows that $E_jm_l=0$ for all $j \in J$
and $1 \leq l \leq n$. By (2) for $\U_J$, and since $m_l \in
L_{\U_J}(\chi)^{\varphi_l}$ for all $l$, the elements $m_l$ are
scalar multiples of $\overline{1} \in L_{\U_J}(\chi)$. Therefore
$m= \Phi(\overline{1} \otimes m'')$, where $m'' \in
L_{\U'}(\chi)^{\chi'}$. Then by Lemma \ref{components} (3),
$0=E_im = \Phi(\overline{1} \otimes E_i m'')$, hence $E_i m'' =0$
for all $i \in \I'$; thus, $m''$ is a scalar multiple of
$\overline{1}$ by (2) for $\U'$. Hence $m\in \ku\overline{1}$ and
$\chi'=\chi$.

\medbreak We next show that (2) for $L(\chi)$ implies (1). Let
$0 \neq M \subset L(\chi)$ be a $\U$-submodule, and let $0 \neq m
\in M$. Then $\U_J\, m$ is a finitely generated $\U_J$-submodule
of $L(\chi)$. By Lemmas \ref{Verma} and \ref{hi} (3), there is $u
\in U_J$ such that $um$ is an element of maximal weight in $\U_J
m$ and $E_j um=0$ for all $j \in J$. Then $\U' um$ is
$\U'$-finitely generated and by the same reason there is an
element $u' \in \U'$ such that $u'um$ is an element of maximal
weight in $\U' um$ and $E_i u'um =0$ for all $i \in \I'$. Then
$u'um$ is a $\ku$-linear combination of elements of the form
$F_{i_1} \cdots F_{i_n}E_{l_1} \cdots E_{l_m}gum$, where $n,m \geq
0$, $i_1,\dots, i_n$, $l_1,\dots,l_m \in \I'$ and $g \in \G$.
Since the elements $F_{i_1},\dots,F_{i_n}$ and $E_{l_1},
\dots,E_{l_m}$ commute or skew-commute with $E_j$ for all $j \in
J$, it follows that $E_j u'um=0$ for all $j \in J$. Let $\chi'$ be
the weight of $u'um$. Since (2) holds for $L(\chi)$ it follows
that $\chi'=\chi$, and $u'um\in \ku \ell_{\chi}$. Hence
$\ell_{\chi} \in M$, and $M =L(\chi)$ since $\ell_{\chi}$
generates the $\U$-module $L(\chi)$.

\medbreak {\rt Finally (3) for $L(\chi)$ follows from the isomorphism $\Phi$ in Lemma \ref{components}.}
\epf

{\rt For an algebra $A$ we denote the set of isomorphism classes of finite-dimensional left $A$-modules by $\Irr(A)$.

\begin{Cor}\label{simplebijection}
\emph{(1)} The map
$$\widehat{\G}^+ \to \{ [L] \mid L \in \mathcal{C}^{hi}, L \text{ integrable  and simple} \},$$
defined by $\chi \mapsto [L(\chi)]$, is bijective.

\emph{(2)} Assume that the Cartan matrix of $\D_{red}$ is of finite type. Then the map in \textup{(1)} defines a bijection
$\widehat{\G}^+ \to \Irr(\U).$
\end{Cor}
\pf (1) The map is well-defined and injective by Corollary \ref{L} and Theorem \ref{mainsimple}. To prove surjectivity let $L \in \mathcal{C}^{hi}$ be integrable and simple. By Lemma \ref{hi}  $L$ has a maximal weight $\chi$, and by Lemma \ref{Verma} and Proposition \ref {dominant} $L \cong L(\chi)$.

(2) By Lemma \ref{components} and the arguments in the proof of Theorem \ref{mainsimple} it suffices to assume that the braiding matrix is of the form $(q^{d_ia_{ij}})$, where $(d_ia_{ij})$ is the symmetrized Cartan matrix of finite type and $0 \ne q \in \ku$ is not a root of unity. Then the claim follows from \cite[5.9, 5.15. 6.26]{J}.
\epf}

\subsection{Complete Reducibility Theorems}\label{sectionreducibil}

\

Here is one of the main results of the present paper {\rt extending \cite[6.2.2]{L}}, the analogue
of (b) in the Introduction.

\begin{Theorem}\label{main}
Let $M$ be an integrable module in $\mathcal{C}^{hi}$. Then $M$ is
completely reducible and $M$ is a direct sum  of $L(\chi)'s$
where $\chi \in \widehat{\G}^+$.
\end{Theorem}
\pf By Theorem \ref{mainsimple} it suffices to show that $M$ is
completely reducible. We proceed by induction on the number of
connected components of $\I$.

\smallbreak Let $\I$ be connected. We may assume that $M \neq 0$.
We need only show that $M$ is a sum of simple $\U$-submodules.
Thus we may further assume that $M$ is $\U$-finitely generated,
and $M=M_C$ for some coset $C$ of $Q$. By Proposition \ref{OG} (1)
the operator $\Omega_G$ for $C$ is locally finite. Since
generalized eigenspaces of module endomorphisms are submodules we
may assume that $M$ is a generalized eigenspace of $\Omega_G$ with
eigenvalue $\lambda$.

\smallbreak Let $N$ be a proper $\U$-submodule of $M$. It suffices
to show that there exists a simple $\U$-submodule $S$ such that $S
\cap N =0$. Then $M$ has a simple submodule (take $N=0$), and $M$
must be the sum of all simple submodules (take $N$ to be this
sum).

\smallbreak Let $m \in M \setminus N$ and set $L=\U \cdot m$. Then
$L/(N \cap L) \neq 0$ is finitely generated and has a maximal
weight $\chi$ by Lemma \ref{hi} (3). Since $\chi$ is also a weight
for $L$ there is a maximal weight $\chi'$ for $L$ which satisfies
$\chi \leq \chi'$. By the characterization of the eigenvalues of
$\Omega_G$ in Proposition \ref{OG} (3) we have that
$G(\chi)=\lambda= G(\chi')$. By Proposition \ref{dominant}  and
Lemma \ref{Verma} both characters $\chi$ and $\chi'$ are dominant.
Hence $\chi=\chi'$ by Lemma \ref{G} (2). Therefore $\chi$ is a
maximal weight for $L$. The projection $L \to L/(N \cap L)$
induces a surjection $L^{\chi} \to (L/(L \cap N))^{\chi}$. We
choose $\ell \in L^{\chi} \setminus N$. Then $S=\U \cdot {\rt \ell}$ is
simple by Lemma \ref{Verma}, Proposition \ref{dominant} and
Theorem \ref{mainsimple}, and $S \cap N =0$.

\smallbreak In the general case let $J$ be a connected component
of $\I$ and $\I'=\I \setminus J$. Let $\U_J$ be the subalgebra of
$\U$ generated by $\G$ and $E_j,F_j$, $j \in J$. Let $\U'$ be the
subalgebra of $\U$ generated by $\G$ and $E_i,F_i$, $i \in \I'$.
We assume that any integrable $\U'$-module in the category
$\mathcal{C}^{hi}$ for $\U'$ is completely reducible.

\smallbreak Let $M$ be a finitely generated and integrable
$\U$-module in $\mathcal{C}^{hi}$, and let $N \subset M$ be a
proper $\U$-submodule. As before it suffices to show that there
exists a simple $\U$-module $S \subset M$ such that $N \cap S =0$.
By the first part of the proof $M$ is completely reducible over
$\U_J$. Hence there exists a simple $U_J$-submodule $S_1 \subset
M$ such that $N \cap S_1 =0$. Let $m \in S_1$ with $S_1 = \U_J m$.
Since $\U'm \not\subset N$ and $\U'm$ is completely reducible by
induction there is a simple $\U'$-module $S_2 \subset \U'm$ such
that $S_2 \cap N =0$. By Theorem \ref{mainsimple} there is a
character $\chi \in \widehat{\G}$ and an element $u \in
\U'^{\chi}$ such that  $S_2= \U' um$ and $E_i um =0$ for all $i
\in \I'$. As in the proof of Theorem \ref{mainsimple} it follows
that $E_j um =0$ for all $j \in J$.  By Proposition \ref{dominant}
and Theorem \ref{mainsimple}  $S:=\U um$ is simple over $\U$.
Moreover $S \cap N =0$, since $S$ is $\U$-simple and $S
\not\subset N$. \epf

\subsection{Reductive pointed Hopf algebras}

\

Let $A$ be an algebra and $B \subset A$ a subalgebra.
We say that

\begin{itemize}
  \item []
$A$ is  {\em reductive} if any
finite-dimensional left $A$-module is completely reducible.

  \item [] $A$
is $B$-{\em reductive} if every finite-dimensional left $A$-module
which is $B$-semisimple when restricted to $B$ is  $A$-semisimple.
\end{itemize}

A pointed Hopf algebra $H$ with $\G=G(H)$ is called $\G$-reductive
if it is $\ku\G$-reductive. Compare with \cite{KSW1, KSW2}.

\begin{Cor}\label{Gammareductive}
$\U$ is $\G$-reductive.
\end{Cor}
\pf This is a special case of Theorem \ref{main} since
finite-dimensional $\U$-modules which are completely reducible
over $k\G$ are  integrable objects of $\mathcal{C}^{hi}$. \epf

\medbreak In Theorem \ref{mainreductive}, we shall need a
generalization of Lemma \ref{EF}. Suppose $q \in \ku^\times$ is not a
root of unity. As usual, we define

$$[a] = \frac{q^a - q^{-a}}{q-q^{-1}},\quad \left[\begin{matrix}
a\\ n \end{matrix}\right]= \frac{[a][a-1]
\cdots[a-n+1]}{[1][2]\cdots[n]}, \quad [n]^{!}=[1][2] \cdots[n]$$

\noindent for all $a,n \in \mathbb{Z}$ and $n >0$, and
$\displaystyle{\left[\begin{matrix}
a\\
0
\end{matrix}\right]=1,\; [0]^{!}=1.}$

\begin{Lem}\label{J}
Let $0 \neq \ell \in \ku$. Let $A$ be an algebra with elements
$E,F,K,L$ such that $K$ and $L$ are invertible and

\begin{align}
KL&=LK,\\
KEK^{-1} &=q^2 E, & KFK^{-1} &=q^{-2} F,\\\label{KF} LEL^{-1}
&=q^2 E, &
LFL^{-1} &=q^{-2} F,\\
EF - FE &= \ell(K - L^{-1}).
\end{align}

Let $(K,L;a) = \displaystyle{\frac{Kq^a -
L^{-1}q^{-a}}{q-q^{-1}}}$, $a \in \mathbb{Z}$. \begin{enumerate}
\item For all $r,s \in \mathbb{Z}, r,s \geq 0,$
\begin{align*}
E^rF^s &= \sum_{i=0}^{\min(r,s)}F^{s-i} h_i(r,s) E^{r-i}, \text{where}\\
h_i(r,s) &={\rt \ell}(q-q^{-1})\left[\begin{matrix}
r\\
i
\end{matrix}\right]
\left[\begin{matrix}
s\\
i
\end{matrix}\right] \left[i\right]^{!}
\prod_{j=1}^i (K,L;i-(r+s)+j).
\end{align*}
\item $KL$ acts semisimply on any finite-dimensional left $A$-module.
\end{enumerate}
\end{Lem}
\pf By rescaling $E$ we may assume that $\ell=(q - q^{-1})^{-1}$.
Let $\lambda \in \ku$ be an eigenvalue of $K$ on $M$. Then by
\eqref{KF} for any natural number $r >0$, $F^r$ maps the
generalized eigenspace of $\lambda$ with respect to $K$ into
$\cup_{n \geq 1} \ker(K - \lambda q^{-2r})^n$. Since $q$ is not a
root of unity there is an integer $s >0$ such that $F^sM=0$.  The
elements $(K,L;a)$ satisfy the same rules as $[K;a]=[K,K;a]$ in
\cite{J}. One can check that the proofs of (1) in \cite[Lemma
1.7]{J} and of (2) in \cite[Prop. 2.3]{J}-- which uses (1)-- work
in our more general situation. In particular,
$\Big(\prod_{j=-(s-1)}^{s-1} (KL - q^{-2j})\Big)M=0.$
\epf

\begin{Rem}\label{standard}
The following is a standard fact in abelian group theory. Let $A$
be a subgroup of an abelian group $B$ with $[B:A] < \infty$. If
$M$ is a ${\rt \ku B}$-module such that $M_{\vert {\rt \ku A}}$ is semisimple then $M$
is semisimple.
\end{Rem}

\pf If $\lambda\in \widehat A$, then we denote by $M_{\lambda}$
the isotypic component of type $\lambda$. Then $M =
\oplus_{\lambda\in \widehat A} M_{\lambda}$ and each $M_{\lambda}$
is a ${\rt \ku B}$-submodule. Thus, we can assume that $M = M_{\lambda}$ for
some $\lambda$. There exists $\Lambda\in \widehat B$ extending
$\lambda$ since the multiplicative group of the algebraically closed field $\ku$ is a divisible hence injective abelian group. Then $A$ acts trivially on $M' =
M\otimes \ku_{\Lambda^{-1}}$, and this becomes a module over the
finite group $B/A$. Hence $M'$ is a semisimple ${\rt \ku B}$-module, and so
is $M \simeq M'\otimes \ku_{\Lambda}$. \epf

Recall that $\G^2$ denotes the subgroup of $\G$ generated by the products $K_1L_1,\dots,K_{\theta}L_{\theta}$, see Definition \ref{regular}.

\begin{Theorem}\label{mainreductive}
{\rt T}he following are equivalent:
\begin{enumerate}
\item ${\bf U}$ is reductive.
\item $[\G : \G^2]$ is finite.
\end{enumerate}
If ${\bf U}$ is reductive, then the Cartan matrix of $\D_{red}$ is invertible.
\end{Theorem}
\pf Suppose $\bf U$ is reductive. There is a well-defined
surjective algebra map $\U \to k[\G/\G^2]$ mapping all $E_i$
and $F_i$ onto zero and any $g \in \G$ onto its residue class in
$\G/\G^2$. Hence the group algebra $k[\G/\G^2]$ is reductive, and
$\G/\G^2$ must be finite.

Conversely suppose that $\G/\G^2$ is finite. Let $M$ be a
finite-dimensional left $\bf U$-module. Then for any $1 \leq i
\leq \theta$, the elements $E_i,F_i,K_i,L_i$ in $\bf U$ satisfy
the assumptions of Lemma \ref{J}. Hence $K_iL_i$ acts semisimply
on $M$ by Lemma \ref{J}. Then we obtain from Remark \ref{standard}
that $\G$ acts semisimply on $M$. Thus $M$ is a semisimple $\bf
U$-module by Corollary \ref{Gammareductive}.

Finally, if (1) and (2) hold, then $(a_{ij})$ is invertible by Lemma \ref{independent} (3).
\epf

{\rt \section{A characterization of quantum groups}\label{characterization}}

We now turn to the representation theory of {\rt the} more general pointed
Hopf algebra $\Uc(\D, \lambda)$ {\rt where $\D$ is generic and of finite Cartan type. Let $\D'$ be the datum with
perfect linking parameter $\lambda'$ associated to the set $\I^s$
of all non-linked vertices as in Theorem \ref{Levi}. Then $ \Uc(\D',\lambda') \cong \Uc(\D_{red},\ell) :=\U$, where $(\D_{red}$, $\ell)$ is} the reduced datum
and its linking parameter associated to $(\D',\lambda')$ as in Lemma
\ref{perfectreduced}. Thus $\D_{red}$ is generic of finite Cartan type, hence regular by Lemma \ref{independent}. By Section \ref{section:Levi} there is a projection of Hopf algebras  $$\pi_{\D} : \Uc(\D,\lambda) \to \U.$$
We may consider then any
{\rt $\Uc(\D_{red},\ell)$-module} as a $\U$-module via $\pi_{D}$, and $\pi_{\D}$ induces a mapping $\pi_{\D}^*$ from the isomorphism classes of $\U$-modules to the isomorphism classes of  $\Uc(\D,\lambda)$-modules.
Let $\G^2 \subset \G$
{\rt be the subgroup defined in Definition \ref{regular} for $\D_{red}$}.

\medskip

\begin{Prop}\label{prop:general-simple}
Let $\D$ be a generic datum of finite Cartan type with linking parameter $\lambda$ and define $\pi_{\D}$ as above. Then
$$ \widehat{\G}^+ \to \Irr((\Uc(\D,\lambda)),\; \chi \mapsto \pi_{\D}^*[L(\chi)],$$
is bijective.
\end{Prop}

\pf This follows from Corollary \ref{simplebijection} and
\cite[Theorem 4.6]{RS2}. \epf

 {\rt \begin{Lem}\label{Lemmaperfect} Let $\D$ be a generic YD-datum
of finite Cartan type and abelian group $\G$, and let $\lambda$ be
a linking parameter for $\D$. Let $h \in \I$, and assume that $h$ is not linked.
Define  $\D',\lambda'$,
$\Uc=\Uc(\D, \lambda)$, $\Uc'=\Uc(\D',\lambda')$  and $K = \Uc^{co\Uc'}$ as in Theorem \ref{Levi} for $L = \{h\}$. Then $M =
\Uc/(\Uc\Uc'^+ + \Uc(K^+)^2)$ is a finite-dimensional vector space.
\end{Lem}}
\pf
{\rt Let $J$ be the connected component of $\I$ containing $h$, and let $X_J = \oplus_{i \in J} \ku x_i$. Then the natural algebra map $\rho : \NA(X) \to \Uc$ is injective by
Theorem \ref{Masuoka}. We view $\rho$ as an inclusion. By Theorem \ref{Levi} $K$ is contained in $\NA(X_J)$. Since $K = \Uc^{co\Uc'}$ is a left coideal subalgebra of $\Uc$, it follows that $K$} {\rt is a left coideal subalgebra of $\NA(X_J)\ku\G \cong \NA(X_J) \# \ku\G$. Hence $K \subset \NA(X_J)$ is a left coideal subalgebra in the braided sense, that is, $\Delta_{\NA(X_J)} (K) \subset \NA(X_J) \otimes K$. Moreover, $K$ is $\mathbb{N}^J$-graded by Theorem \ref{Levi}. Since the braiding of $X_J$ is of finite Cartan type, it follows from \cite[Corollary 6.16]{HS} that there are finitely many $\mathbb{N}^J$-homogeneous elements $a_1,\dots a_m \in K$ such that for all $i \in J$ the subalgebra $\ku(a_i)$ is isomorphic to $\NA(\ku a_i)$, and the multiplication map $\ku(a_m) \otimes \cdots \ku(a_1) \to K$ is bijective. Since $\NA( X_J)$ is an integral domain (see for example \cite[Theorem 4.3]{AS-crelle}), for each i the Nichols algebra $\NA(\ku a_i)$ is a polynomial ring. Hence the elements $a_m^{n_m}\cdots a_1^{n_1}, n_1,\dots,n_m \geq0$, form a $\ku$-basis of $K$. The existence of such a PBW-basis can also be derived from \cite{Kh}. Since $M$ is an epimorphic image of $K/(K^+)^2$ by the decomposition $K \# \Uc' \cong \Uc$ in Theorem \ref{Levi}, it follows that $M$ is the $\ku$-span of the images of $a_1,\dots,a_m$ thus finite-dimensional.}
\epf

\begin{Theorem}\label{pointedreductive}
Let $\D$ be a generic YD-datum of finite Cartan type with abelian
group $\G$, and let $\lambda$ be a linking YD-datum for $\D$.

\emph{(i)} The following are equivalent:
\begin{enumerate}
\item $\Uc(\D,\lambda)$ is $\G$-reductive.
\item The linking parameter $\lambda$ of $\D$ is perfect.
\end{enumerate}

\medbreak
\emph{(ii)} The following are equivalent:
\begin{enumerate}
\item $\Uc(\D,\lambda)$ is reductive.
\item
\begin{enumerate}
\item The linking parameter $\lambda$ of $\D$ is perfect.
\item $[\G : \G^2]$ is finite.
\end{enumerate}
\end{enumerate}

\end{Theorem}
\pf (i) We assume that $U(\D,\lambda)$ is $\G$-reductive, and that the linking is not perfect.
We choose an element $h \in \I$ which is not linked and define $L = \{h \}$. Let  $M =\Uc/(\Uc\Uc'^+ + \Uc(K^+)^2)$ as in Lemma \ref{Lemmaperfect}. {\rt Then $M$ is finite-dimensional by Lemma \ref{Lemmaperfect}.} By Corollary \ref{M} $x_{h}M\neq 0$. Hence $M$ is not semisimple
since by \cite[Theorem 4.6]{RS2} any finite-dimensional simple
${\rt U(\D,\lambda)}$-module is annihilated by $x_{h}$.

To obtain a contradiction we finally show that $M$ is semisimple as a $\G$-module by
restriction. The vector space $\Uc(\D,\lambda)$ is the $\ku$-span
of elements of the form $xh$, $x$ a monomial in the elements
$x_1,\dots,x_{\theta}$, $h \in \G$. Let $g \in \G$, then $gxh =
\chi(g) xhg$ for some $\chi \in \widehat{\G}$. Hence in the module
$M$, we have $g\overline{xh} = \chi(g) \overline{xhg} = \chi(g)
\overline{xh}$ since $g-1 \in {U'}^{+}$. Thus $M$ as a $\G$-module
is the sum of weight spaces.

Conversely
assume that the linking parameter is perfect. Then $U(\D,\lambda)
\cong U(\D_{red},\ell)$ for some generic reduced YD-datum. Since
the Cartan matrices of $\D$ and $\D_{red}$ are of finite type,
$\D_{red}$ is regular by Lemma \ref{independent} (2), and
$U(\D_{red},\ell)$ is $\G$-reductive by Corollary
\ref{Gammareductive}.

(ii) follows from the argument in the proof of (i) and Theorem \ref{mainreductive}.
\epf

Theorem \ref{pointedreductive} combined with \cite[Theorem 1.1]{AA}, that generalizes the main result of \cite{AS-crelle},
gives immediately the following characterization of quantized enveloping algebras.

\begin{Theorem}\label{fingrowth-lifting}
Let $H$ be a pointed Hopf algebra with finitely generated abelian
group $G(H)$, and generic infinitesimal braiding. Then the following are equivalent:
\begin{enumerate}
\item $H$ is a $\G$-reductive domain with finite Gelfand-Kirillov dimension.
\item The group $\Gamma := G(H)$ is free abelian of finite rank,
and there exists a {\rt reduced generic datum of finite Cartan type ${\mathcal D}_{red}$ for $\Gamma$ with linking parameter $\ell$
such that
$H \simeq \U=\Uc({\mathcal D}_{red}, \lambda)$ as Hopf algebras.}
\end{enumerate}
{\rt If $H$ satisfies \textup{(2)}, then $H$ is reductive iff $[\G : \G^2]$ is finite.}
\end{Theorem}


\end{document}